%% file: mainV1.tex
\documentclass[12pt,preprint]{elsarticle}
\usepackage{graphicx} % Required for inserting images
\usepackage{amssymb,amsfonts} % For \mathbb
\usepackage{amsmath} % For pmatrix and math environments
\usepackage[margin=2.5cm]{geometry}

\usepackage{bm}
\usepackage{mathptmx,mathrsfs,mathtools}
\usepackage{stmaryrd}
\usepackage{hyperref}
\usepackage{pifont}
\usepackage{cleveref}
\usepackage{nicematrix}
\usepackage{xcolor}
\usepackage{float}
\usepackage{tcolorbox}
\usepackage{cancel}
\usepackage{tabularx}
\usepackage{array}
\usepackage{subfigure}
\usepackage{subcaption}

\newtheorem{rmk}{Remark}

\usepackage{xcolor}
\usepackage{tcolorbox}
\usepackage{tikz}
\usepackage{todonotes}

\biboptions{sort&compress}

\begin{document}
\begin{frontmatter}

%\title{A unified algebraic treatment of scalar interface problems}

\title{An algebraic reformulation of scalar interface problems}

\author[iitm]{M Arrutselvi\fnref{fnlab2}}
\author[iitm]{Sundararajan Natarajan\corref{cor1}\fnref{fnlab1}}   
\ead{snatarajan@iitm.ac.in}

\address[iitm]{Department of Mechanical Engineering, Indian Institute of Technology Madras, Chennai 600036, India}

\fntext[fnlab2]{NBHM Postdoctoral Fellow}
\cortext[cor1]{Corresponding author}
\fntext[fnlab1]{Pandurangan Faculty Fellow}

\begin{abstract}
This work presents a unified algebraic reformulation of scalar interface problems based on a shared primary field and an independent localized jump field. The two-sided interface traces are represented using a single primary interface variable and an auxiliary variable describing their jump, thereby retaining a single set of geometric interface nodes without conventional node duplication or partition-of-unity enrichment. The resulting augmented block-matrix system explicitly organizes the bulk, interface, and coupling contributions. The formulation accommodates continuous and discontinuous interface behaviour, including prescribed jumps and linear or nonlinear constitutive transmission conditions. For linear bulk problems with nonlinear interface laws, the nonlinear residual and consistent tangent contributions are confined to the interface-associated blocks, leaving the bulk stiffness contributions unchanged. An appropriate diagonal row-column scaling is also investigated to improve the conditioning of the system in the stiff-interface limit. Numerical studies using analytical and manufactured solutions assess the convergence and performance of the formulation. The proposed framework provides a discretization-independent algebraic representation for incorporating diverse scalar interface conditions within a common computational structure.
\end{abstract}


\end{frontmatter}
 
%---
\section{Introduction}
\input{Intro2}

%--- model problem, governing equation and weak form
\input{modelProb}
\input{unifiedReform}

% \section{Unified algebraic treatment}
% \label{sec:newformu}

% Run a quick condition number comparison between your $\mathbf{\delta }$-formulation and standard double-node FEM across $\beta \in [10^{-6}, 10^{6}]$.

% \begin{equation} 
% \begin{aligned}
% u^{(2)} &= \sum_{a \in \Omega_2\setminus \Gamma} N_a u^{(2)}_a+\sum_{b \in \Gamma} N_b u^{(2)}_b, \qquad x \in \Omega_2\\
% u^{(1)} &= \sum_{\alpha \in \Omega_1\setminus\Gamma} N_\alpha u^{(1)}_\alpha + \sum_{b \in \Gamma} N_b u_b^{(2)} - \sum_{b\in \Gamma} \overline{N}_{b}\delta_b, \qquad x \in \Omega_1
% \end{aligned}
% \end{equation} 

\section{Numerical Examples}
\label{sec:numexam} 
The numerical performance of the proposed algebraic reformulation is assessed through a series of examples involving diffusion across internal interfaces subject to the different transmission conditions listed in \Cref{tab:interface_conditions}. Particular attention is given to the representation of interface discontinuities through auxiliary nodal jump variables and the reuse of bulk finite element matrices when changing the interface conditions. The accuracy and convergence properties of the formulation are investigated through one- and two-dimensional examples.
%---------------
\subsection{Errors}
 Let each subdomain $\Omega_i$ is discretized independently using a mesh $\mathcal{T}_h^{(i)}$, such that
\begin{equation*}
\overline{\Omega}_i
= \cup_{K\in\mathcal{T}_h^{(i)}}\overline{K},
\qquad i=1,2,
\end{equation*}
where $h_i=\max\limits_{K\in\mathcal{T}_h^{(i)}}h_K$ denotes the mesh size in $\Omega_i$, and $h=\max\{h_1,h_2\}$ is the global mesh parameter. Let $u^{(i)}$ and $u_h^{(i)}$ denote the exact and numerical solutions in $\Omega_i$, respectively. The corresponding approximation errors are defined as $e^{(i)}=u^{(i)}-u_h^{(i)}$.

The numerical errors are quantified using the broken $L^2$-norm and $H^1$-seminorm, defined by
\begin{equation*}
\begin{aligned}
\|e\|_{L^2(\Omega)}^2
&=
\sum_{i=1}^{2}
\sum_{K\in\mathcal{T}_h^{(i)}}
\int_K
(e^{(i)})^\top e^{(i)}\,\mathrm{d}\Omega,\\
|e|_{H^1(\Omega_1\cup\Omega_2)}^2
&=
\sum_{i=1}^{2}
\sum_{K\in\mathcal{T}_h^{(i)}}
\int_K
(\nabla e^{(i)})^\top
\nabla e^{(i)}\,\mathrm{d}\Omega.
\end{aligned}
\end{equation*}
These norms measure the approximation errors in the bulk subdomains independently and remain well-defined even when the numerical solution exhibits a discontinuity across $\Gamma$.
\subsection{One-dimensional heat conduction across an imperfect interface}

We consider a one-dimensional steady-state heat conduction problem in a composite domain consisting of two homogeneous materials separated by an internal interface. This example illustrates the incorporation of different transmission conditions within the proposed algebraic formulation and assesses the resulting numerical solutions against their analytical counterparts. The governing equation is
\begin{equation}
-\frac{\mathrm{d}}{\mathrm{d}x}
\left(k(x)\frac{\mathrm{d}u}{\mathrm{d}x}\right)=0,
\qquad x\in(0,1)\setminus\{x_b\},
\end{equation}
where the thermal conductivity is piecewise constant:
\begin{equation}
k(x)=
\begin{cases}
k^{(1)}, & x\in\Omega_1=(0,x_b),\\
k^{(2)}, & x\in\Omega_2=(x_b,1).
\end{cases}
\end{equation}
Here, $\Gamma=\{x_b\}$ denotes the internal interface separating the two materials. The prescribed boundary temperatures are
\begin{equation}
u(0)=u_L,\qquad u(1)=u_R.
\label{eqn:onedBC}
\end{equation}

The heat flux is defined by
\begin{equation}
q^{(i)}=-k^{(i)}\frac{\mathrm{d}u^{(i)}}{\mathrm{d}x},
\qquad i=1,2,
\end{equation}
with the jump convention
\begin{equation}
\llbracket u\rrbracket
=u^{(2)}(x_b)-u^{(1)}(x_b),
\qquad
\llbracket q\rrbracket
=q^{(2)}(x_b)-q^{(1)}(x_b).
\end{equation}
Five interface conditions are considered: (i) perfect thermal contact, with continuous temperature and heat flux; (ii) continuous temperature with a prescribed heat-flux jump; (iii) a prescribed temperature jump with continuous heat flux; (iv) prescribed jumps in both temperature and heat flux; and (v) an imperfect thermal interface governed by a Kapitza-type constitutive law. The corresponding analytical solutions are derived in \Cref{appx:oned_anal} and serve as reference solutions for assessing the numerical results.

For the spatial discretization, each subdomain is represented by one linear finite element, with the interface located at the common element boundary. The nodal coordinates are
\[
x_0=0,\qquad x_1=x_b,\qquad x_2=1.
\]
Following the proposed algebraic representation, we retain a single primary interface unknown $u_b$ and introduce an independent jump variable
\begin{equation}
u_b=u_b^{(2)},\qquad
\delta_b=u_b^{(2)}-u_b^{(1)},
\end{equation}
so that
\begin{equation}
u_b^{(1)}=u_b-\delta_b.
\end{equation}
The global vector of unknowns is therefore arranged as
\begin{equation}
\mathbf{U}=
\begin{bmatrix}
u_0 & u_b & u_2 & \delta_b
\end{bmatrix}^{\top}.
\end{equation}
Thus, the two physical interface traces are represented without introducing geometrically duplicated interface nodes.

For linear interpolation, the element stiffness matrices are
\begin{equation}
\mathbf{K}^{(1)}
=c_1
\begin{bmatrix}
1&-1\\
-1&1
\end{bmatrix},
\qquad
\mathbf{K}^{(2)}
=c_2
\begin{bmatrix}
1&-1\\
-1&1
\end{bmatrix},
\end{equation}
where
\begin{equation}
c_1=\frac{k^{(1)}}{x_b},
\qquad
c_2=\frac{k^{(2)}}{1-x_b}.
\end{equation}
For a Kapitza-type interface with resistance parameter $\beta$, the constitutive law is
\begin{equation}
\delta_b=-\beta q,
\end{equation}
where $q$ denotes the common heat flux across the interface. The corresponding interface contribution is characterized by $\rho=1/\beta$. Combining the two bulk element contributions with the interface contribution gives the augmented system
\begin{equation}
\begin{bNiceMatrix}
c_1 & -c_1 & 0 &
\Block[draw=green,fill=green!8]{3-1}{
\begin{gathered}
c_1\\[-2pt]
-c_1\\[-2pt]
0
\end{gathered}
}\\
-c_1 & c_1+c_2 & -c_2 & \\
0 & -c_2 & c_2 & \\
\Block[draw=blue,fill=blue!8]{1-3}{
\begin{gathered}
c_1 \quad -c_1 \quad 0
\end{gathered}
}
&&&
\Block[draw=red,fill=red!8]{1-1}{
c_1+\dfrac{1}{\beta}
}
\end{bNiceMatrix}
\begin{bmatrix}
u_0\\
u_b\\
u_2\\
\delta_b
\end{bmatrix}
=
\begin{bmatrix}
f_0\\
f_b\\
f_2\\
0
\end{bmatrix}.
\label{eqn:finalstiff}
\end{equation}
Here, $f_0$, $f_b$, and $f_2$ denote the corresponding nodal load contributions. The prescribed heat-flux jump, when present, contributes an additional term $g_q$ to the right-hand side associated with the primary interface equation.

\begin{rmk}
The augmented matrix in \Cref{eqn:finalstiff} explicitly identifies the coupling between the bulk degrees of freedom and the independent interface jump variable. The upper-left $3\times3$ block corresponds to the assembled bulk finite element stiffness matrix, whereas the last row and column represent the coupling associated with $\delta_b$. The last diagonal entry contains the contribution $c_1$ from the bulk-jump coupling and the interface conductance contribution $\rho=1/\beta$. Thus, the interface constitutive response is incorporated through the additional jump variable without modifying the underlying bulk element matrices.

For prescribed interface jumps, the corresponding constraints can be imposed directly by algebraic elimination, while a prescribed flux jump is incorporated into the primary interface equation. Consequently, the same bulk element contributions can be retained when switching between different interface conditions, with only the interface-related algebraic equations and right-hand-side contributions being modified.
\end{rmk}

Enforcing the Dirichlet boundary conditions in \Cref{eqn:onedBC}, the reduced system for the interface unknowns becomes
\begin{equation}
\begin{bmatrix}
c_1+c_2 & -c_1\\
-c_1 & c_1+\rho
\end{bmatrix}
\begin{bmatrix}
u_b\\
\delta_b
\end{bmatrix}
=
\begin{bmatrix}
c_1u_L+c_2u_R+g_q\\
-c_1u_L
\end{bmatrix},
\label{eqn:stiff_reduced}
\end{equation}
where $g_q=0$ in the absence of a prescribed heat-flux jump. The five interface conditions are now examined individually.

\paragraph{Case 1: Perfect interface ($\delta_b=0$, $g_q=0$)}
For a perfect interface, both the temperature and heat flux are continuous. Imposing $\delta_b=0$ in \Cref{eqn:stiff_reduced} gives
\begin{equation}
(c_1+c_2)u_b=c_1u_L+c_2u_R.\nonumber
\end{equation}
Consequently, the interface temperature is
\begin{equation}
u_b=\frac{c_1u_L+c_2u_R}{c_1+c_2}.
\label{eq:perfect_ub}
\end{equation}
The heat fluxes on the two sides are
\begin{equation}
q^{(1)}=c_1(u_L-u_b),
\qquad
q^{(2)}=c_2(u_b-u_R).\nonumber
\end{equation}
Substitution of \Cref{eq:perfect_ub} yields
$\llbracket q\rrbracket=q^{(2)}-q^{(1)}=0$,
thereby recovering the perfect-interface conditions.

\paragraph{Case 2: Prescribed heat-flux jump ($\delta_b=0$, $g_q$ prescribed)}
In this case, the temperature remains continuous, while a prescribed jump in the heat flux is introduced. Setting $\delta_b=0$ in \Cref{eqn:stiff_reduced} gives
\begin{equation}
(c_1+c_2)u_b
=c_1u_L+c_2u_R+g_q.\nonumber
\end{equation}
Hence,
\begin{equation}
u_b=
\frac{c_1u_L+c_2u_R+g_q}{c_1+c_2}.\label{interf2u2}
\end{equation}
Using \Cref{interf2u2} in the corresponding heat fluxes,
\begin{equation*}
q^{(1)}=c_1(u_L-u_b),
\qquad
q^{(2)}=c_2(u_b-u_R),
\end{equation*}
 gives~$ 
\llbracket q\rrbracket=g_q$.~Thus, a prescribed interface source or sink produces a discontinuity in the heat flux while maintaining temperature continuity.

\paragraph{Case 3: Prescribed temperature jump ($\delta_b=g_u$, $g_q=0$)}
Here, a prescribed temperature jump is imposed, while the heat flux remains continuous. Substituting $\delta_b=g_u$ and $g_q=0$ into \Cref{eqn:stiff_reduced} gives
\begin{equation}
(c_1+c_2)u_b-c_1g_u
=c_1u_L+c_2u_R.\nonumber
\end{equation}
Therefore,
\begin{equation}
u_b=
\frac{c_1u_L+c_2u_R+c_1g_u}{c_1+c_2}.\label{interf3}
\end{equation}
Then substituting \Cref{interf3} in the corresponding heat fluxes
\begin{equation}
q^{(1)}=c_1(u_L-u_b+g_u),
\qquad
q^{(2)}=c_2(u_b-u_R),\label{flux-case3}
\end{equation}
satisfy~$\llbracket q\rrbracket=0.$~This case represents a prescribed strong discontinuity in temperature accompanied by continuous heat flux.

\paragraph{Case 4: Simultaneous prescribed temperature and heat-flux jumps ($\delta_b=g_u$, $g_q$ prescribed)}
In this case, both interface jumps are prescribed. Substituting $\delta_b=g_u$ into \Cref{eqn:stiff_reduced} gives
\begin{equation}
(c_1+c_2)u_b-c_1g_u
=c_1u_L+c_2u_R+g_q.\nonumber
\end{equation}
The primary interface value is consequently
\begin{equation}
u_b=
\frac{c_1u_L+c_2u_R+g_q+c_1g_u}
{c_1+c_2}.\label{interf4}
\end{equation}
The corresponding heat fluxes are again given by \Cref{flux-case3}. Now, putting \Cref{interf4} into \Cref{flux-case3}  satisfy~$ \llbracket q\rrbracket=g_q.$~Thus, the formulation accommodates simultaneous prescribed discontinuities in the temperature and heat flux.

\paragraph{Case 5: Imperfect interface (Kapitza condition)}
We now consider an imperfect thermal interface governed by the constitutive relation
\begin{equation}
\delta_b=-\beta q,
\qquad
\llbracket q\rrbracket=0,
\end{equation}
where $\beta>0$ denotes the interface thermal resistance. In this case, the temperature jump is an unknown, determined by the interface constitutive law. Setting $g_q=0$ and $\rho=1/\beta$, the reduced system in \Cref{eqn:stiff_reduced} directly provides the primary interface temperature and the jump variable. Defining
\begin{equation}
\mathcal{D}_{\beta}
=c_1c_2+\rho(c_1+c_2),\nonumber
\end{equation}
the solution is
\begin{equation}
 u_b=
\frac{c_1c_2u_R+\rho(c_1u_L+c_2u_R)}
{\mathcal{D}_{\beta}}\qquad\text{and}\qquad
\delta_b=
\frac{c_1c_2(u_R-u_L)}
{\mathcal{D}_{\beta}}.\nonumber
% \label{eq:kapitza_solution}
\end{equation}
The corresponding heat flux is continuous and satisfies the prescribed Kapitza relation. In particular, as $\beta\to0$, the interface approaches the perfect-contact limit, with $\delta_b\to0$. Thus, the perfect and imperfect interface behaviours are accommodated within the same algebraic representation, with the former recovered as a limiting case of the latter.

These five cases demonstrate the applicability of the proposed formulation to different transmission conditions within a common algebraic structure. The computed solutions can be verified against the corresponding analytical solutions derived in \Cref{appx:oned_anal}, with which they agree. Importantly, the bulk element stiffness contributions remain unchanged across all five cases. The different interface behaviours are incorporated through appropriate algebraic constraints, constitutive relations, and modifications to the interface-associated equations, illustrating the flexibility of the proposed reformulation.

%%%%---- one d nonlinear Kapitza
\subsection{One-dimensional nonlinear Kapitza interface problem}

The proposed algebraic formulation can also accommodate nonlinear interface
constitutive laws without modifying the underlying bulk finite element
discretization. The nonlinear behaviour enters through the interface residual
and its consistent tangent, while the bulk stiffness contributions and
bulk--interface coupling structure remain unchanged; see
\Cref{eq:nonlinear_jacobian_system}. Thus, the nonlinear extension requires
updating only the interface-associated contributions within an iterative
solution procedure. Unlike the preceding examples, for which analytical
solutions were available, the present problem is solved numerically. A
manufactured solution is employed to assess the accuracy of the formulation
in representing nonlinear temperature jumps across an internal interface.

Consider the domain \(\Omega=(0,1)\), partitioned into
\(\Omega_1=(0,x_b)\) and \(\Omega_2=(x_b,1)\), with
\(\Gamma=\{x_b\}\) and \(x_b=1/2\). The governing problem is
\begin{equation}
\begin{aligned}
-\frac{\mathrm{d}}{\mathrm{d}x}
\left(k^{(i)}\frac{\mathrm{d}u^{(i)}}{\mathrm{d}x}\right)
&=0
&&\text{in }\Omega_i,\quad i=1,2,\\
\llbracket q\rrbracket&=0
&&\text{at }x=x_b,\\
\llbracket u\rrbracket
&=-\bigl(1+u_b^2\bigr)q
&&\text{at }x=x_b,\\
u^{(1)}(0)&=\frac54,\qquad
u^{(2)}(1)=-\frac32,
\end{aligned}
\label{eq:nonlinear_kapitza_problem}
\end{equation}
where \(k^{(1)}=2\), \(k^{(2)}=1\),
\(q^{(i)}=-k^{(i)}\,\mathrm{d}u^{(i)}/\mathrm{d}x\), and
\(q=q^{(1)}=q^{(2)}\) denotes the continuous interface flux, oriented in
the positive \(x\)-direction. The primary interface trace and the
temperature jump are defined as
\[
u_b=u^{(2)}(x_b),\qquad
\llbracket u\rrbracket
=u^{(2)}(x_b)-u^{(1)}(x_b).
\]

The manufactured exact solution is
\begin{equation}
u_{\mathrm{ex}}(x)=
\begin{cases}
\dfrac{5-2x}{4},
& 0\leq x<\dfrac12,\\[1em]
-\dfrac{1+2x}{2},
& \dfrac12<x\leq 1.
\end{cases}
\label{eq:nonlinear_kapitza_exact}
\end{equation}
It gives \(u^{(1)}(x_b)=1\), \(u_b=u^{(2)}(x_b)=-1\), and
\(\llbracket u\rrbracket=-2\). Moreover, the flux is continuous with
\(q=1\), and the nonlinear interface law is satisfied exactly.

For the numerical approximation, each subdomain is discretized using one
two-noded linear element, with nodes at \(x=0\), \(x=x_b\), and \(x=1\).
The global unknown vector is
\[
\mathbf z=
\begin{Bmatrix}
u_0&u_b&u_2&\delta_\Gamma
\end{Bmatrix}^{\top},
\]
where \(u_b=u_b^{(2)}\) and
\(\delta_\Gamma=u^{(2)}(x_b)-u^{(1)}(x_b)\).

\noindent Using the nonlinear residual formulation in
\Cref{eq:nonlinear_residual_system}, the assembled residual is
\begin{equation}
\mathbf R(\mathbf z)=
\begin{Bmatrix}
4u_0-4u_b+4\delta_\Gamma\\
-4u_0+6u_b-2u_2-4\delta_\Gamma\\
-2u_b+2u_2\\
4u_0-4u_b+
\bigl(4+M(u_b)\bigr)\delta_\Gamma
\end{Bmatrix},\nonumber
% \label{eq:nonlinear_kapitza_residual}
\end{equation}
where ~$ M(u_b)=\dfrac{1}{1+u_b^2}.$  

The first three components represent the bulk equilibrium equations,
while the fourth enforces the nonlinear interface constitutive condition.
The corresponding consistent Jacobian is
\begin{equation}
\mathbf J(\mathbf z)=
\begin{bmatrix}
4&-4&0&4\\
-4&6&-2&-4\\
0&-2&2&0\\
4&-4+\delta_\Gamma M'(u_b)&0&4+M(u_b)
\end{bmatrix},
\nonumber%\label{eq:nonlinear_kapitza_jacobian}
\end{equation}
where $ ~ M'(u_b)=-\dfrac{2u_b}{(1+u_b^2)^2}.$  

After imposing the prescribed boundary values
\(u_0=5/4\) and \(u_2=-3/2\), perform the Newton-Raphson method's $n^{th}$  iteration solve of the system, 
\begin{equation}
\mathbf J(^{(n)}\mathbf z)
\Delta\,^{(n)}\mathbf z
=-\mathbf R(^{(n)}\mathbf z),
\qquad
  ^{(n+1)}\mathbf z
= ^{(n)}\mathbf z+\Delta ^{(n)}\mathbf z.\nonumber
%\label{eq:nonlinear_kapitza_newton}
\end{equation}

\begin{table}[htbp]
\centering
\renewcommand{\arraystretch}{1.05}
\setlength{\tabcolsep}{12pt}
\caption{Newton convergence history for the one-dimensional nonlinear
Kapitza interface problem.}
\label{tab:nonlinearkapitzaresults}
\begin{tabular}{c|c|c|c}
\hline
Iteration & $u_1$ & $\delta_\Gamma$ &
$\|\mathbf R\|_2$\\
\hline
1 & $-0.71428571$ & $-1.57142857$ & $5.000\times10^{0}$\\
2 & $-0.97854740$ & $-1.9678211$ & $5.309\times10^{-1}$\\
3 & $-0.97854740$ & $-1.96782111$ & $3.766\times10^{-2}$\\
4 & $-0.99993038$ & $-1.99989557$ & $1.218\times10^{-4}$\\
5 & $-1.00000000$ & $-2.00000000$ & $1.212\times10^{-9}$\\
6 & $-1.00000000$ & $-2.00000000$ & $1.776\times10^{-15}$\\
\hline
\end{tabular}
\end{table}

Starting from \(^{(0)}u_b=0\) and
\(^{(0)}\delta_\Gamma=0\), the Newton iteration converges to
\(u_b=-1\) and \(\delta_\Gamma=-2\) within five iterations.
As shown in \Cref{tab:nonlinearkapitzaresults}, the reduced residual norm
decreases to approximately \(1.78\times10^{-15}\), demonstrating convergence
to the manufactured solution.

The nonlinear interface behaviour is incorporated through a localized
residual and its consistent tangent, while the bulk stiffness contributions
and coupling structure remain unchanged from the linear Kapitza problem.
Furthermore, the original three-node mesh is retained: the interface
discontinuity is represented by the auxiliary algebraic unknown
\(\delta_\Gamma\), without geometrically duplicating the interface node.
Thus, the proposed formulation accommodates the nonlinear interface law
without introducing separate geometric interface nodes for the two
subdomains.

%-----2D example
\subsection{Two-dimensional heat conduction across a straight interface}

\begin{figure}[H]
\centering
\subfigure[]{
    \includegraphics[scale=0.5]{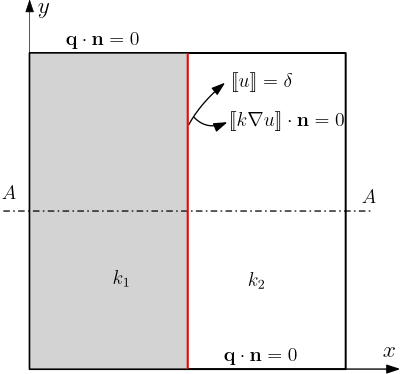}
    \label{fig:planargeom}
}
\hspace{50pt}
\subfigure[]{
    \includegraphics[scale=0.5]{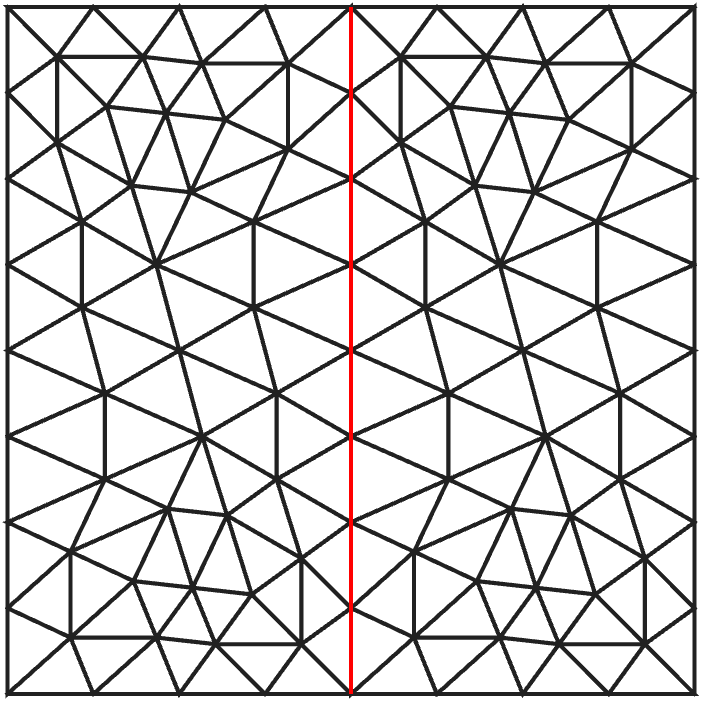}
    \label{fig:planarfemesh}
}
\caption{Domain with a straight internal interface:
(a) problem geometry and boundary conditions and
(b) representative computational mesh. The red line denotes the interface.}
\label{fig:2DprobDescription}
\end{figure}

We next consider heat conduction in a two-dimensional domain containing
an internal straight interface. Three test cases are investigated:
(a) a constant prescribed temperature jump (Case A),
(b) thermal contact resistance (Kapitza interface, Case B), and
(c) a spatially varying prescribed temperature jump (Case C).
The computational domain is the unit square
\[
\Omega=(0,1)\times(0,1),
\]
partitioned into two non-overlapping subdomains,
\[
\Omega_1=(0,x_b)\times(0,1),\qquad
\Omega_2=(x_b,1)\times(0,1),
\]
with the interface
\[
\Gamma=\{(x,y)\in\Omega:x=x_b\},\qquad x_b=0.5.
\]
The thermal conductivities in the two subdomains are denoted by
\(k_1\) and \(k_2\), respectively. Throughout this section, the
temperature jump and the normal heat-flux jump are defined as
\[
\llbracket u\rrbracket_\Gamma
=u^{(2)}-u^{(1)},\qquad
\llbracket q_n\rrbracket_\Gamma
=q_n^{(2)}-q_n^{(1)},
\]
where
\[
q_n^{(i)}=-k^{(i)}\nabla u^{(i)}\cdot\mathbf n
\]
and \(\mathbf n\) is the unit normal directed from \(\Omega_1\)
towards \(\Omega_2\). Figure~\ref{fig:2DprobDescription} shows the
problem geometry and a representative computational mesh.

\paragraph{Case A: Heat conduction with a prescribed temperature jump}

In this case, homogeneous Neumann boundary conditions are imposed on
the bottom and top boundaries, \(y=0\) and \(y=1\), respectively.
A constant temperature jump is prescribed across the interface, while
the normal heat flux remains continuous:
\begin{equation}
\llbracket u\rrbracket_\Gamma=\delta_\Gamma,
\qquad
\llbracket q_n\rrbracket_\Gamma=0.
\label{eqn:Strint}
\end{equation}
Here, \(\delta_\Gamma\) denotes the prescribed temperature jump.
The thermal conductivities are chosen such that
\[
\frac{k^{(2)}}{k^{(1)}}=6,
\]
and the source term is \(f=1\). The corresponding analytical solution,
which is independent of \(y\), is given by~\cite{hansbohansbo2002}:
\begin{equation}
u_{\mathrm{ex}}(x,y)=
\begin{cases}
\displaystyle
\frac{(3k^{(1)}+k^{(2)})x}{4{(k^{(1)})}^2+4k^{(1)}k^{(2)}}
-\frac{x^2}{2k^{(1)}},
& x<x_b,\\[2ex]
\displaystyle
\frac{(k^{(2)}-k^{(1)})+(3k^{(1)}+k^{(2)})x}
     {4(k^{(2)})^2+4k^{(1)}k^{(2)}}
-\frac{x^2}{2k^{(2)}}
+\delta_\Gamma,
& x>x_b.
\end{cases}
\label{eq:caseA_exact}
\end{equation}
The Dirichlet boundary values on \(x=0\) and \(x=1\) are prescribed
consistently with this analytical solution.

\begin{figure}[htpb]
\centering
\subfigure[]{
    \includegraphics[scale=0.5]{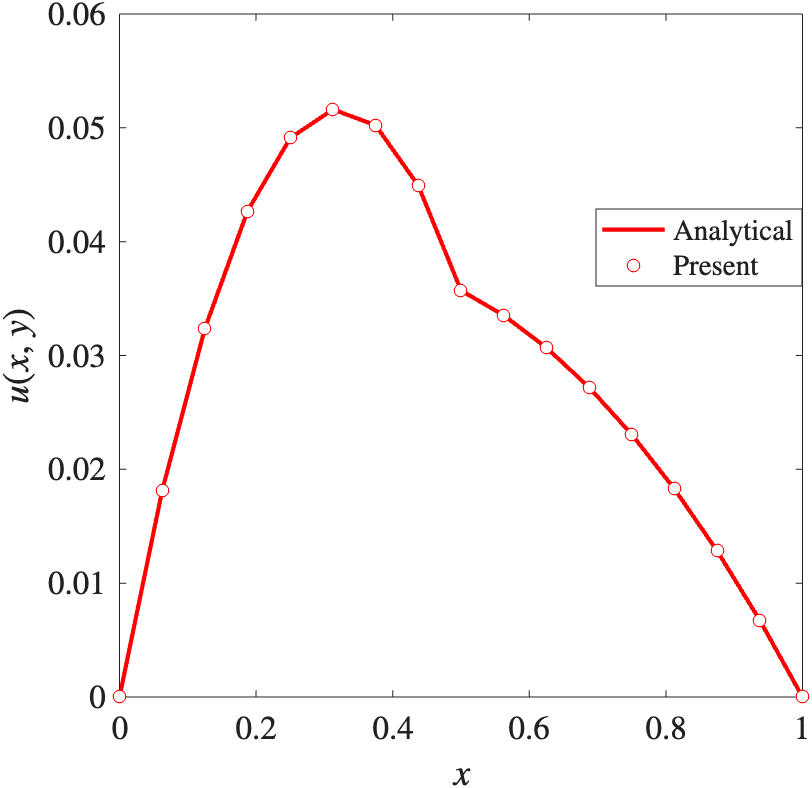}
    \label{fig:lpltnojump}
}
\hspace{50pt}
\subfigure[]{
    \includegraphics[scale=0.5]{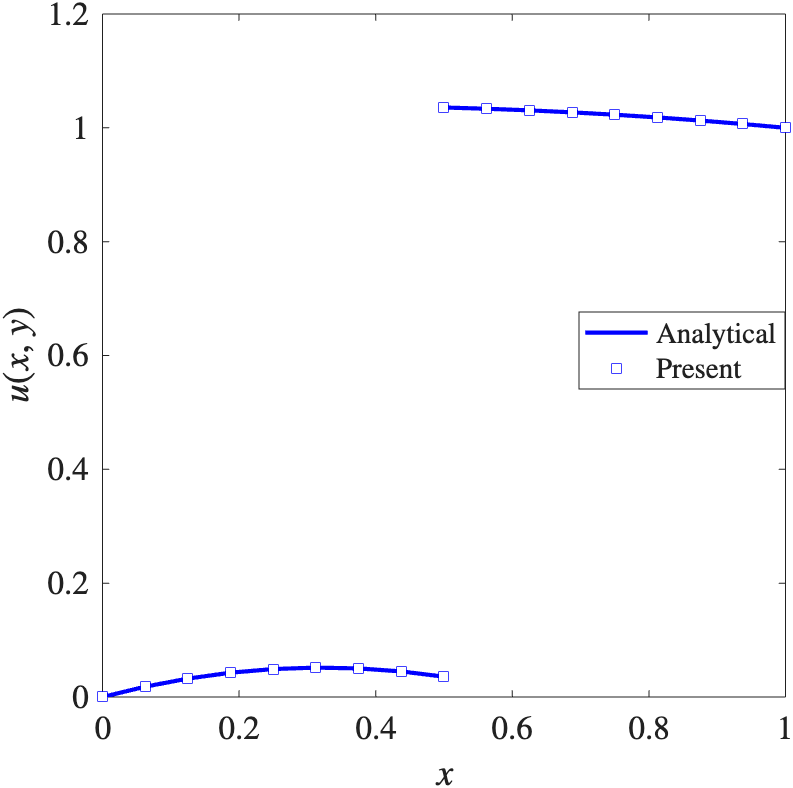}
    \label{fig:lpltjump}
}
\subfigure[]{
    \includegraphics[scale=0.5]{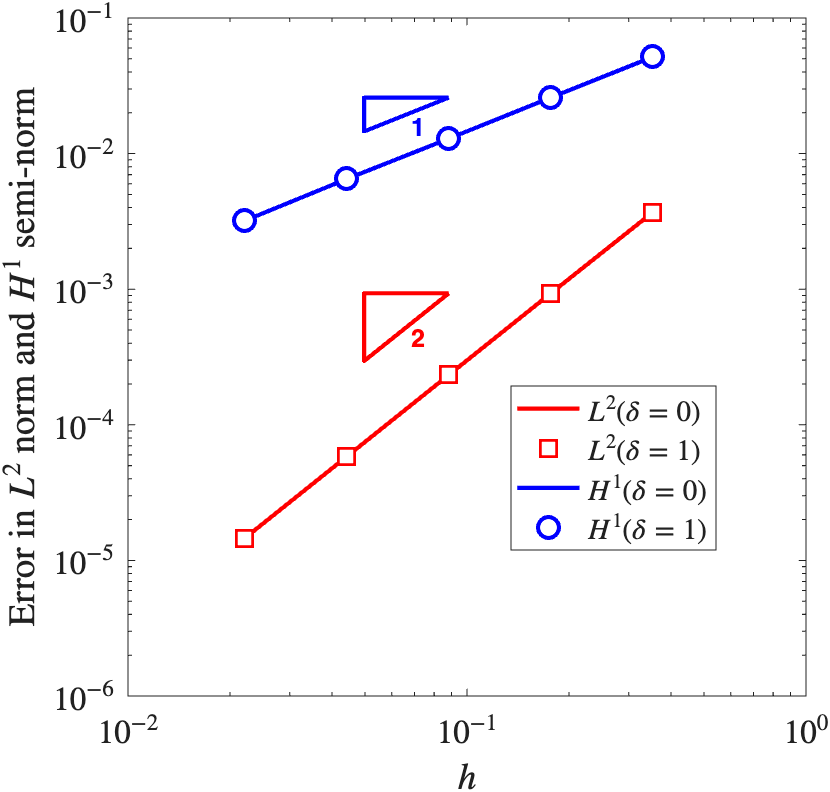}
    \label{fig:2dl2h1}
}
\caption{Case A: Two-dimensional heat conduction with a prescribed
temperature jump. (a,b) Variation of \(u(x,y=1/2)\) for
\(\delta_\Gamma=0\) and \(\delta_\Gamma=1\), respectively.
(c) Convergence of the errors in the \(L^2\) norm and the
\(H^1\) seminorm for the two prescribed-jump cases.}
\label{fig:twod_jump}
\end{figure}

Figure~\ref{fig:twod_jump}(a,b) shows the variation of
\(u(x,y=1/2)\) for \(\delta_\Gamma=0\) and
\(\delta_\Gamma=1\), respectively. The numerical results obtained
using the proposed algebraic reformulation agree well with the
analytical solution. The convergence of the errors in the \(L^2\)
norm and the \(H^1\) seminorm under mesh refinement is presented in
Figure~\ref{fig:twod_jump}(c). The observed convergence rates are
two in the \(L^2\) norm and one in the \(H^1\) seminorm, consistent
with the expected optimal rates for the approximation.

\paragraph{Case B: Heat conduction with Kapitza resistance}

We next consider heat conduction across an interface with thermal
contact resistance. As in Case A, homogeneous Neumann boundary
conditions are imposed on the top and bottom boundaries. The
temperature is prescribed on the left and right boundaries as
\[
u^{(1)}(0,y)=T_L=0,\qquad
u^{(2)}(1,y)=T_R=1.
\]
The interface condition is the Kapitza law
\begin{equation}
\llbracket u\rrbracket_\Gamma=-\beta q_n,
\label{eq:caseB_kapitza}
\end{equation}
where \(\beta\) is the thermal contact resistance and \(q_n\)
denotes the common normal heat flux. For the analytical solution,
it is convenient to introduce the positive-gradient flux
\[
p_s=k^{(i)}\frac{\partial u^{(i)}}{\partial x}=-q_n.
\]
The interface condition in terms of this quantity is therefore
\(\llbracket u\rrbracket_\Gamma=\beta p_s\).
The analytical solution is
\begin{equation}
u_{\mathrm{ex}}(x,y)=
\begin{cases}
\displaystyle
\frac{p_s}{k^{(1)}}(x-x_0)+T_L,
& x<x_b,\\[1.5ex]
\displaystyle
\frac{p_s}{k^{(2)}}(x-x_1)+T_R,
& x>x_b,
\end{cases}
\label{eq:caseB_exact}
\end{equation}
where
\[
p_s=
\frac{T_R-T_L}
{\displaystyle
\frac{x_b-x_0}{k^{(1)}}
+\frac{x_1-x_b}{k^{(2)}}
+\beta},
\qquad
x_0=0,\quad x_1=1.
\]

\begin{figure}[htpb]
\centering
\includegraphics[scale=0.5]{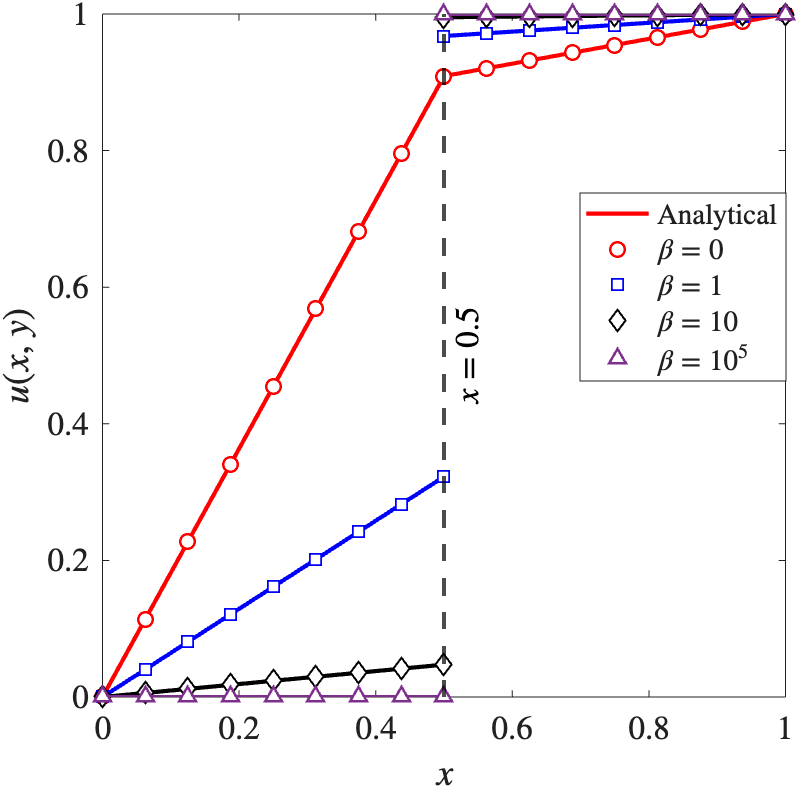}
\caption{Case B: Variation of \(u(x,y=1/2)\) for different values
of the thermal contact resistance \(\beta\).}
\label{fig:twodcaseb}
\end{figure}

Figure~\ref{fig:twodcaseb} shows the variation of
\(u(x,y=1/2)\) for different values of the thermal contact
resistance. The numerical results agree well with the analytical
solution.

The errors in the \(L^2\) norm and the \(H^1\) seminorm for different
values of \(\beta\), as well as for different mesh sizes with
\(\beta=1\), are reported in Table~\ref{tab:combined_errors}.
For all the considered resistance values and mesh resolutions,
the errors remain close to machine precision. This behaviour is
consistent with the fact that the exact solution is piecewise
linear and is therefore represented exactly by the linear
approximation on the considered meshes, up to numerical
round-off errors.

\begin{table}[htbp]
\centering
\caption{Case B: Errors in the \(L^2\) norm and the \(H^1\)
seminorm for (a) different thermal contact resistances with
\(h=0.1768\), and (b) different mesh sizes with \(\beta=1\).}
\label{tab:combined_errors}

\begin{minipage}{0.48\textwidth}
\centering
\begin{tabular}{lcc}
\hline
\(\beta\) & \(L^2\) error & \(H^1\) error \\
\hline
\(1\times10^{-8}\) & \(1.9446\times10^{-16}\) & \(9.7140\times10^{-16}\) \\
\(1\)              & \(1.6759\times10^{-16}\) & \(8.2918\times10^{-16}\) \\
\(10\)             & \(1.7257\times10^{-16}\) & \(7.2463\times10^{-16}\) \\
\(1\times10^{5}\)  & \(2.8684\times10^{-16}\) & \(1.1161\times10^{-15}\) \\
\hline
\end{tabular}
\caption*{(a) Variation with thermal contact resistance.}
\label{tab:kapitza_errors}
\end{minipage}
\hfill
\begin{minipage}{0.48\textwidth}
\centering
\begin{tabular}{ccc}
\hline
\(h\) & \(L^2\) error & \(H^1\) error \\
\hline
0.3536 & \(1.7257\times10^{-16}\) & \(7.2463\times10^{-16}\) \\
0.1768 & \(9.2826\times10^{-16}\) & \(3.1486\times10^{-15}\) \\
0.0884 & \(3.7586\times10^{-15}\) & \(1.3275\times10^{-14}\) \\
0.0442 & \(7.7640\times10^{-15}\) & \(2.9284\times10^{-14}\) \\
\hline
\end{tabular}
\caption*{(b) Variation with mesh size.}
\label{tab:mesh_errors}
\end{minipage}
\end{table}

\begin{figure}[htpb]
\centering
\subfigure[Double-node formulation with a Lagrange multiplier]{
    \includegraphics[scale=0.5]{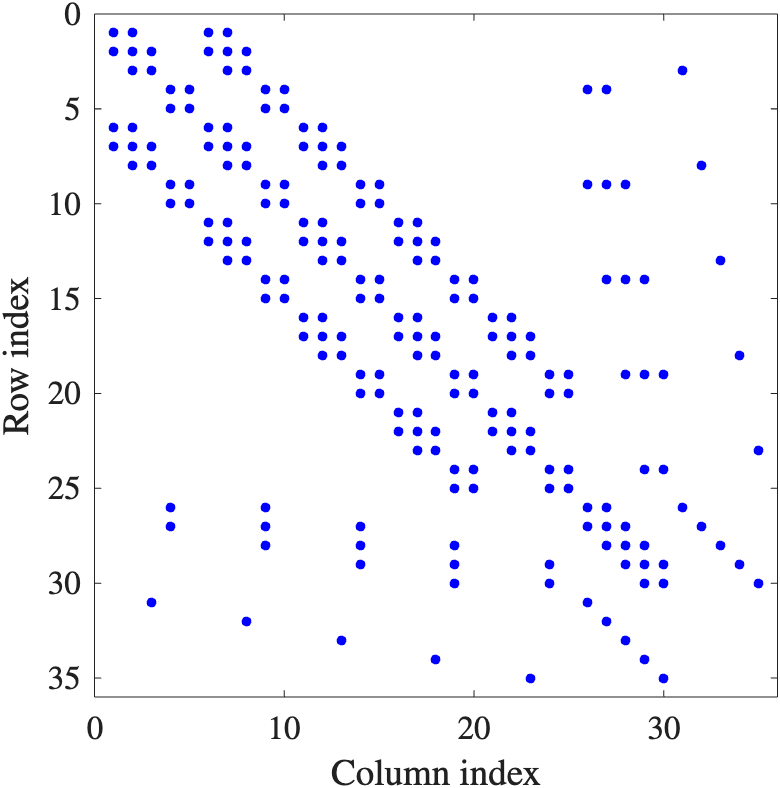}
}
\hspace{50pt}
\subfigure[Proposed algebraic reformulation]{
    \includegraphics[scale=0.5]{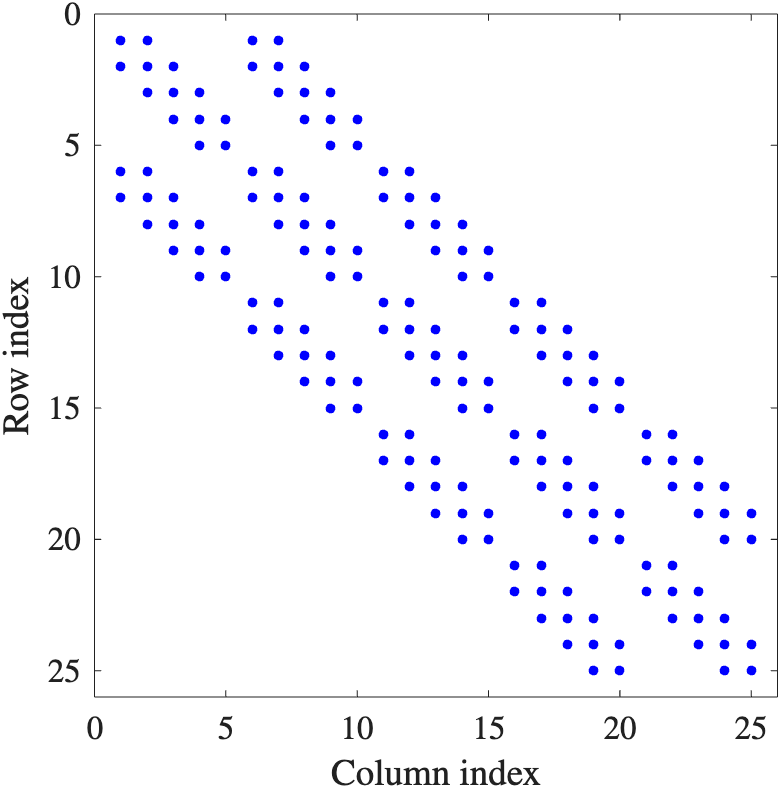}
}
\caption{Comparison of the nonzero entries in the system matrices
for the double-node formulation with a Lagrange multiplier and
the proposed algebraic reformulation. The respective numbers of
nonzero entries are 202 and 169.}
\label{fig:nnz_stiff}
\end{figure}

For comparison, the above cases are also solved using a conventional
double-node formulation in which the interface nodes are duplicated
and the interface constraint is enforced using a Lagrange multiplier.
Both formulations produce the same discrete solution and convergence
behaviour. For the \(4\times4\) quadrilateral mesh shown in
Figure~\ref{fig:nnz_stiff}, the proposed reformulation retains the
banded structure of the bulk stiffness matrix, whereas the
Lagrange-multiplier formulation introduces a saddle-point block.
In the double-node formulation, each interface node is associated
with three unknowns: two physical interface values and one
Lagrange multiplier. In the proposed formulation, these are
replaced by two unknowns: one shared interface value and one
independent jump variable. Thus, for a prescribed interface jump,
the jump variable can be eliminated directly, without introducing
a Lagrange-multiplier unknown or the associated saddle-point
structure.

We emphasize that the proposed formulation is not obtained through
a Schur complement reduction of the Lagrange-multiplier system.
Rather, it is an algebraic reparameterization of the two interface
traces: the duplicated interface values are represented by a
shared interface value and an explicit jump variable. When the
jump is prescribed, the jump variable is eliminated directly,
thereby avoiding the additional Lagrange-multiplier unknown and
the associated saddle-point structure.

\paragraph{Case C: Heat conduction with a spatially varying prescribed jump}

Finally, we consider a spatially varying prescribed temperature jump
along the interface. The thermal conductivities are taken to be
\(k^{(1)}=k^{(2)}=1\), and the governing equation is
\begin{equation}
-\nabla\cdot(k^{(i)}\nabla u^{(i)})=f^{(i)}
\qquad\text{in }\Omega_i,\quad i=1,2.
\label{eq:caseC_governing}
\end{equation}
A manufactured solution is prescribed as
\begin{equation}
u_{\mathrm{ex}}(x,y)=
\begin{cases}
u^{(1)}(x,y)=x,
& (x,y)\in\overline{\Omega}_1,\\[1ex]
u^{(2)}(x,y)=x+\sin^2(\pi y),
& (x,y)\in\overline{\Omega}_2.
\end{cases}
\label{eq:caseC_exact}
\end{equation}
The corresponding temperature jump is
\begin{equation}
\llbracket u\rrbracket_\Gamma
=g_u(y),\qquad
g_u(y)=\sin^2(\pi y).
\label{eq:caseC_jump}
\end{equation}
Substitution of the exact solution into
\eqref{eq:caseC_governing} gives the manufactured source term
\begin{equation}
f^{(i)}(x,y)=
\begin{cases}
0, & (x,y)\in\Omega_1,\\[1ex]
-2\pi^2\cos(2\pi y), & (x,y)\in\Omega_2.
\end{cases}
\label{eq:caseC_source}
\end{equation}
The Dirichlet boundary conditions on the left and right boundaries
are prescribed from the exact solution:
\[
u^{(1)}(0,y)=0,\qquad
u^{(2)}(1,y)=1+\sin^2(\pi y).
\]
Homogeneous Neumann boundary conditions are imposed on the top and
bottom boundaries.

\begin{figure}[htpb]
\centering
\subfigure[]{
    \includegraphics[scale=0.5]{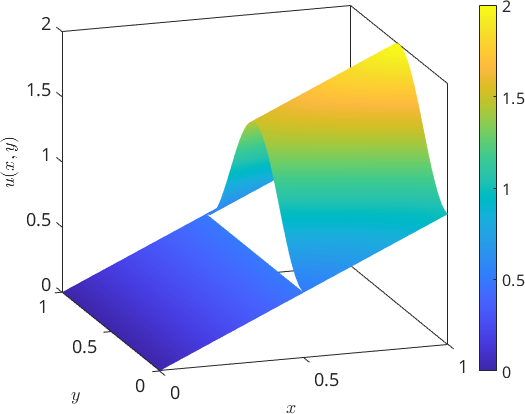}
    \label{fig:casec_cont}
}
\hspace{25pt}
\subfigure[]{
    \includegraphics[scale=0.5]{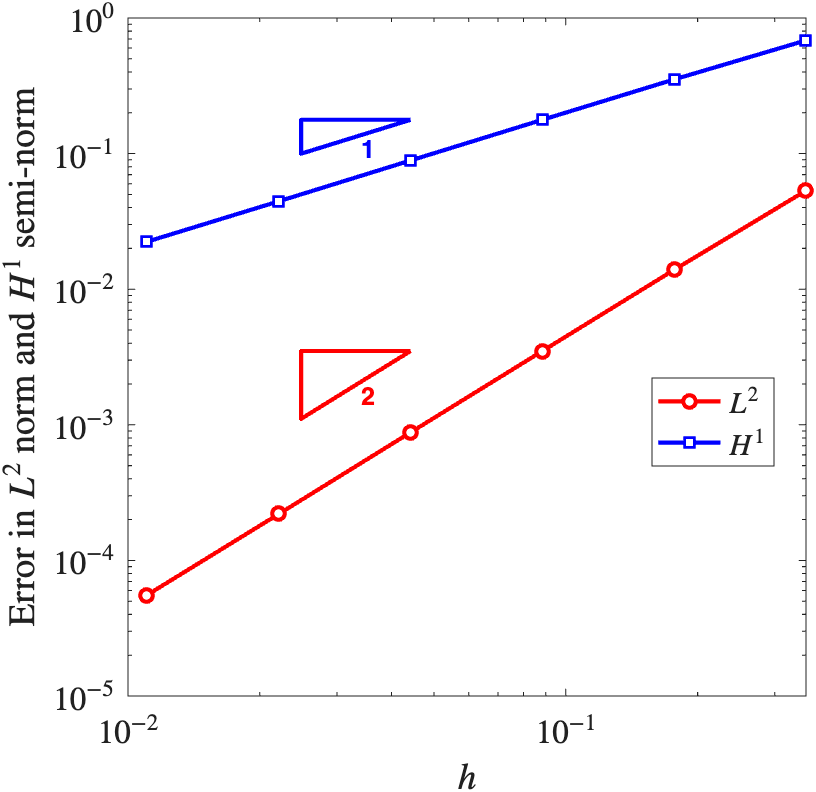}
    \label{fig:casec_conv}
}
\caption{Case C: Two-dimensional heat conduction with a spatially
varying prescribed temperature jump:
(a) contour plot of \(u(x,y)\) and
(b) convergence of the errors in the \(L^2\) norm and the
\(H^1\) seminorm under mesh refinement.}
\label{fig:casec_result}
\end{figure}

Figure~\ref{fig:casec_cont} shows the contour plot of the computed
temperature field, illustrating the spatial variation of the
solution jump along the interface. The convergence of the errors
in the \(L^2\) norm and the \(H^1\) seminorm under mesh refinement
is presented in Figure~\ref{fig:casec_conv}. The observed
convergence rates are two in the \(L^2\) norm and one in the
\(H^1\) seminorm, consistent with the expected optimal rates.
This example demonstrates that the proposed formulation can
accommodate a non-constant prescribed jump along the interface.
In particular, the spatially varying jump is imposed through the
local jump degrees of freedom, while the bulk discretization
remains unchanged.

%--------------------------
\subsection{Conditioning analysis for a two-dimensional Kapitza interface problem}

We investigate the conditioning of the proposed jump-parameterized
formulation for a two-dimensional diffusion problem in
\(\Omega=(0,1)\times(0,1)\), partitioned into
\[
\Omega_1=(0,x_b)\times(0,1),\qquad
\Omega_2=(x_b,1)\times(0,1),
\qquad x_b=\frac12.
\]
The diffusion coefficients are \(k^{(1)}=1\) and \(k^{(2)}=6\).
Across the interface \(\Gamma\), the normal flux is continuous and
the temperature jump satisfies the Kapitza law
\[
q_n^{(1)}=q_n^{(2)}=q_n,\qquad
u_\Gamma^{(1)}-u_\Gamma^{(2)}=\beta q_n,
\]
where \(q_n^{(i)}=-k^{(i)}\nabla u^{(i)}\cdot\mathbf n\) and
\(\beta\) is the interface resistance. The limits
\(\beta\to0\) and \(\beta\to\infty\) correspond to perfect contact
and an insulating interface, respectively.

The domain is discretized using a \(4\times4\) mesh of bilinear
quadrilateral elements. We compare the spectral condition numbers
of the conventional double-node formulation and the proposed
algebraic reformulation for
\[
\bar\beta=\frac{\beta k_{\mathrm{ref}}}{h_\Gamma}
\in\{10^{-6},10^{-3},10^{-1},1,10^3,10^6\},
\]
where \(h_\Gamma\) is a representative interface mesh size and
\(k_{\mathrm{ref}}\) is a reference diffusion coefficient.

As reported in Table~\ref{tab:comparisoncond}, the unscaled
formulation exhibits a deterioration in conditioning as the
interface resistance approaches zero, owing to the dominant
\(1/\beta\) contribution in the interface block. To address this,
we introduce a diagonal scaling of the jump degrees of freedom.
Specifically, we define
\[
\widehat{\boldsymbol{\delta}}_\Gamma
=s\boldsymbol{\delta}_\Gamma,\qquad
s=\sqrt{\frac{\bar\beta}{1+\bar\beta}},
\]
and obtain the scaled stiffness matrix through the congruence
transformation
\begin{equation}
\widehat{\mathbf K}=\mathbf S^\top\mathbf K\mathbf S,
\label{eq:scaled_stiffness}
\end{equation}
where
\[
S_{ii}=
\begin{cases}
1, & i\text{ corresponds to a bulk degree of freedom},\\
s, & i\text{ corresponds to a jump degree of freedom}.
\end{cases}
\]
This change of variables does not alter the physical solution,
which is recovered using
\(\boldsymbol{\delta}_\Gamma
=s^{-1}\widehat{\boldsymbol{\delta}}_\Gamma\).

\begin{figure}[htpb]
\centering
\includegraphics[scale=0.6]{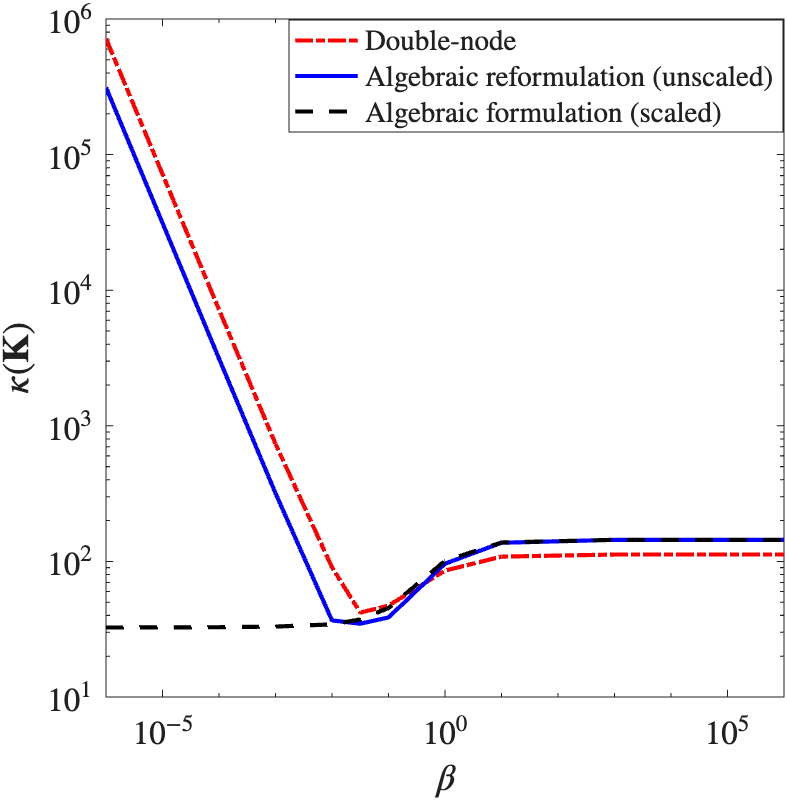}
\caption{Condition numbers as functions of the Kapitza resistance
\(\beta\) for the double-node formulation and the proposed
algebraic reformulation, with and without diagonal scaling.}
\label{fig:condnumbervari}
\end{figure}

The resulting condition numbers are listed in
Table~\ref{tab:comparisoncond} and plotted in
Figure~\ref{fig:condnumbervari}. For small interface resistance,
the scaling substantially reduces the condition number. In
particular, for \(\beta=10^{-6}\), it decreases from
\(3.13\times10^5\) to \(32.56\), compared with
\(7.17\times10^5\) for the double-node formulation. Over the
considered parameter range, the scaled condition number remains
comparable to the bulk-dominated conditioning level.

\begin{table}[htbp]
\renewcommand{\arraystretch}{1.3}
\setlength{\tabcolsep}{10pt}
\centering
\caption{Spectral condition numbers of the conventional double-node
formulation and the proposed algebraic reformulation, with and
without diagonal scaling, for different interface resistance values.}
\label{tab:comparisoncond}
\begin{tabular}{lccc}
\hline
$\beta$ & Double-node formulation
& \multicolumn{2}{c}{Algebraic reformulation}\\
\cline{3-4}
& & Unscaled & Scaled\\
\hline
$10^{-6}$ & $7.1743\times10^5$ & $3.1289\times10^5$ & $32.5561$\\
$10^{-3}$ & $735.8553$ & $318.2447$ & $33.0657$\\
$10^{-1}$ & $47.1389$ & $38.5944$ & $45.3342$\\
$1$       & $85.5444$ & $96.1933$ & $100.8541$\\
$10$      & $108.4306$ & $137.1960$ & $137.8588$\\
$10^3$    & $112.3622$ & $144.2371$ & $144.2441$\\
$10^6$    & $112.4038$ & $144.3122$ & $144.3122$\\
\hline
\end{tabular}
\end{table}

The effectiveness of this scaling stems from the explicit jump
variable in the proposed formulation, whose contribution is
localized to the interface-associated degrees of freedom. Thus,
the jump variables can be scaled without directly rescaling the
bulk degrees of freedom. In contrast, the duplicated interface
unknowns in the double-node formulation also participate in the
bulk stiffness contributions, making such localized scaling less
direct. Although other preconditioning strategies may be developed
for the double-node formulation, the proposed parameterization
provides a simple means of balancing the interface and bulk
contributions, particularly in the small-resistance regime.

\section{Conclusions}

This work introduces a unified algebraic reformulation for the numerical treatment of scalar interface problems involving both continuous and discontinuous primary fields. The central idea is to represent the interfacial discontinuity through an explicit jump variable, while retaining a shared bulk interface field and the original mesh topology. Separating the bulk solution from its interfacial jump provides a common representation for prescribed-jump conditions and constitutively governed interfaces without introducing problem-specific enrichment functions. Compared to conventional double-node formulations and their associated Lagrange multiplier treatments, the proposed approach incorporates prescribed jumps directly through the jump degrees of freedom, avoiding additional multiplier unknowns and the resulting saddle-point structure.

Numerical experiments involving prescribed constant jumps, spatially varying jumps, and thermal contact resistance demonstrate the ability of the formulation to reproduce analytical solutions across distinct interfacial regimes. The observed optimal convergence rates in the \(L^2\) norm and \(H^1\) seminorm confirm the accuracy of the discretization, while comparisons with a standard double-node approach highlight its computational advantages in interface degree-of-freedom count and algebraic system structure. Importantly, the proposed formulation acts as an algebraic reparameterization of the interface degrees of freedom rather than a Schur complement reduction or a specific choice of bulk discretization.

Although linear finite elements are used in the present numerical investigation, the reformulation is independent of the polynomial order of the bulk approximation. Higher-order bulk discretizations can be accommodated by introducing a compatible approximation for the interface jump without altering the underlying bulk approximation space. This feature also provides scope for extending the framework beyond conventional finite elements to other discretization paradigms, including immersed and virtual element methods.

Several extensions merit further investigation. These include three-dimensional interface problems, moving interfaces and Stefan-type phase-change problems, and coupled multiphysics applications such as thermoelastic delamination, cohesive interfaces with thermal contact resistance, and hydraulic fracturing. The treatment of non-conforming discretizations through mortar-type coupling and the investigation of parallel scalability using domain-decomposition methods, such as FETI, are also of interest. Such developments would extend the applicability of the proposed reformulation to more complex interface-dominated simulations while preserving its defining feature: the explicit algebraic representation of interfacial discontinuities through auxiliary jump variables, without modifying the bulk approximation space or duplicating interface nodes.

\section*{Declarations}
\subsection*{Funding}
\noindent Not applicable.

\subsection*{Acknowledgments}
\noindent M.~Arrutselvi is thankful for financial support from the National
Board for Higher Mathematics (NBHM), Department of Atomic Energy,
Government of India (Grant 0204/8/2025/R\&D-II/9073).\\
AI tools may have served as writing and editing assistants; however,
the authors carefully reviewed the manuscript and take full
responsibility for the content of this work.

\subsection*{Conflicts of interest/Competing interests}
\noindent The authors declare no conflict of interests.

\subsection*{Ethics approval and consent to participate}
\noindent Not applicable. This study does not involve human participants, animals, or any personal data.

\subsection*{Availability of data and materials}
\noindent No datasets were generated or analyzed during the current study. The numerical examples presented in this paper were generated as part of the study and are reproducible from the methods and algorithms described. The implementation code is available from the corresponding author upon reasonable request.

\subsection*{Consent for publication}
\noindent Not applicable.

%-- references
\bibliographystyle{elsarticle-num-names}
\bibliography{ref}

% \clearpage 
% \newpage 

\appendix
\input{onedAnal}

\end{document}

%% file: Intro2.tex
Interface problems governed by partial differential equations with discontinuous coefficients, singular sources, or non-standard transmission conditions arise in a wide range of engineering and physical applications. Examples include imperfect thermal and electrical contact in heterogeneous materials, diffusion across material interfaces, transport through semi-permeable membranes, and electrostatic interfaces. Depending on the underlying interface physics, the primary scalar field \(u\) and its normal flux \(q_n\) may remain continuous, exhibit prescribed discontinuities, or satisfy constitutive relations connecting the interfacial flux to the jump in the primary field. Consequently, these problems encompass both weak and strong discontinuities and may involve linear or nonlinear transmission conditions. Their numerical treatment therefore requires a suitable representation of the two-sided physical fields that can accommodate different interface behaviours without compromising computational efficiency.

Several numerical strategies have been developed to represent strong discontinuities in computational discretizations. Conventional double-node techniques and zero-thickness interface or cohesive elements introduce independent degrees of freedom on either side of the interface by duplicating interface nodes. These approaches have been widely used in fracture, delamination, and contact modelling \cite{needleman1987continuum,ortiz1999finite}. Although straightforward, node duplication modifies the mesh connectivity and associated data structures. Moreover, in penalty-based implementations, the near-perfect interface limit may require increasingly large penalty or conductance parameters, potentially resulting in ill-conditioned algebraic systems. An alternative is provided by partition-of-unity enrichment techniques, including the extended finite element method (XFEM) \cite{moes1999finite,belytschko2001arbitrary,friesbelytschko2006,chinsukumar2019,caorialmoreno2024} and the generalized finite element method (GFEM) \cite{babuvska1997partition,strouboulis2000design}. These methods enrich the approximation space with additional functions to represent discontinuities independently of element boundaries. However, they require additional approximation functions and degrees of freedom, together with corresponding modifications to numerical integration and assembly procedures. Other approaches, including phantom-node and embedded formulations \cite{hansbo2004finite}, provide alternative representations without conventional geometric node duplication. Collectively, these developments demonstrate the importance of the algebraic representation of interface unknowns in the numerical treatment of discontinuous fields.

In parallel, general constitutive and thermodynamic frameworks have been developed to describe imperfect-interface behaviour. In particular, Javili et al. \cite{javilikaessmair2014} presented a thermodynamically consistent formulation for general imperfect thermal interfaces, encompassing different interface responses through constitutive relations involving discontinuities in temperature and normal heat flux. Related studies have investigated general strong and weak formulations for imperfect interfaces in coupled multifield problems and their numerical implementation using extended finite element techniques \cite{guliu2014,liugu2014}. These contributions provide comprehensive physical and variational descriptions of imperfect interfaces. The present work addresses a complementary computational question: 
\begin{center} 
``{\it whether the two-sided scalar fields can be represented through an alternative algebraic parameterization that provides a systematic organization of the discrete unknowns and their associated operators}".
\end{center}

The central idea of this work is to introduce an independent localized jump variable alongside a single primary interface field. Consider an interface \(\Gamma\) separating two subdomains, \(\Omega_1\) and \(\Omega_2\). Rather than treating the solution traces \(u^{(1)}\) and \(u^{(2)}\) as two independently introduced interface nodal fields, we retain the interface unknown associated with \(\Omega_2\) as the primary variable and define an auxiliary variable representing the jump between the two traces:
\[
u_b^{(2)}=u_b,\qquad
\delta_b=u_b^{(2)}-u_b^{(1)}.
\]
The trace on the first side is consequently reconstructed as
\[
u_b^{(1)}=u_b-\delta_b.
\]
Thus, the two physical interface traces are represented through a single primary interface field and an independent localized jump field. The jump is therefore an explicit algebraic unknown rather than a quantity evaluated by subtracting two separately approximated physical fields. This provides a direct representation of the interfacial discontinuity, while retaining the original geometric interface-node arrangement. The additional algebraic degrees of freedom arise from the jump variables and not from geometric node duplication.

This representation leads to an augmented global algebraic system whose unknown vector is organized as
\[
\mathbf U=
\begin{bmatrix}
\mathbf u_I^{(1)}\\
\mathbf u_I^{(2)}\\
\mathbf u_b\\
\boldsymbol{\delta}_b
\end{bmatrix},
\]
where \(\mathbf u_I^{(1)}\) and \(\mathbf u_I^{(2)}\) denote the interior degrees of freedom in the two subdomains, respectively, \(\mathbf u_b\) represents the primary interface unknowns, and \(\boldsymbol{\delta}_b\) contains the corresponding nodal jump variables. Substitution of this representation into the discrete governing equations produces a block-matrix system that explicitly identifies the bulk, primary interface, jump, and coupling contributions. Different transmission conditions can then be incorporated through appropriate interface equations and algebraic constraints, without changing the underlying bulk approximation structure. In particular, continuous-field behaviour is recovered by setting \(\boldsymbol{\delta}_b=\mathbf 0\), prescribed jumps can be imposed directly, and unknown jumps governed by linear or nonlinear constitutive interface laws can be treated through the corresponding interface equations.

An additional computational advantage arises when the bulk equations are linear but the interface response is nonlinear. In this setting, the nonlinear residual and consistent tangent contributions associated with the interface constitutive law can be assembled within the interface-related blocks of the global system, leaving the underlying bulk stiffness contributions unchanged. The resulting organization separates the nonlinear interfacial response from the bulk discretization and provides a convenient basis for the numerical treatment of different interface models. Furthermore, since the jump is represented explicitly, its numerical evaluation does not rely on subtracting two nearly equal approximations of the primary field, which can be advantageous when the interfacial discontinuity is small.

The proposed approach is distinct from conventional double-node techniques in that the two physical traces are represented by a primary interface variable and an independent jump variable rather than by geometrically duplicated nodes. It also differs from partition-of-unity enrichment methods, since the discontinuity is introduced through localized algebraic unknowns instead of additional parent-element enrichment functions. More recently, Chou \cite{chou2026} proposed a scalar interface reduction for nonlinear interface problems by decomposing the solution into a continuous component and a global unit-jump response mode, thereby reducing the nonlinear interface condition to a scalar equation for a single jump amplitude. In contrast, the present formulation introduces independent jump unknowns at the interface nodes, allowing the discontinuity to vary spatially along the interface. It therefore does not require the construction of a global unit-jump response mode and can accommodate spatially varying prescribed or constitutively determined jumps within the same augmented algebraic system.

Since the proposed construction is based on the algebraic relationship between the two interface traces rather than on a particular bulk approximation, it is not intrinsically restricted to a specific discretization technique. Subject to the availability of compatible bulk and interface approximation operators, the framework can be incorporated into different numerical methods, including finite element, extended finite element, and virtual element discretizations. The algebraic system also provides a setting for investigating conditioning in the stiff-interface limit. In particular, an appropriate diagonal row-column scaling is considered to improve the numerical behaviour of the system when the bulk and interfacial stiffnesses exhibit substantially different magnitudes.

The main objective of this work is to develop and assess a unified algebraic framework for scalar interface problems through an alternative representation of the two-sided interface unknowns. The principal contributions are summarized as follows:
\begin{itemize} 
\item[\ding{113}] \textit{Alternative algebraic representation}:\, A single primary interface field and an independent localized jump field are introduced to represent the two physical solution traces, avoiding geometric duplication of interface nodes while retaining the algebraic degrees of freedom required to describe discontinuities.

\item[\ding{113}] \textit{Unified block-matrix formulation}:\, The discrete governing equations are organized into an augmented block system that explicitly identifies the bulk, interface, jump, and coupling contributions. Different transmission conditions are incorporated within this common algebraic structure.

\item[\ding{113}] \textit{Direct treatment of interface discontinuities}:\, The jump is treated as an independent unknown, providing a direct representation of prescribed and unknown discontinuities. Continuous-field behaviour is recovered as a special case.

\item[\ding{113}] \textit{Localized treatment of nonlinear interface laws}:\, For linear bulk subdomains coupled through nonlinear interface constitutive relations, the nonlinear residual and consistent tangent contributions are confined to the interface-associated blocks, leaving the bulk stiffness contributions unchanged.

\item[\ding{113}] \textit{Numerical assessment and conditioning}:\, The formulation is examined through representative interface problems, convergence studies, nonlinear interface simulations, and conditioning analysis in the stiff-interface limit. The effectiveness of an appropriate diagonal row-column scaling is also investigated.
\end{itemize}

The remainder of the paper is organized as follows. Section 2 introduces the scalar interface problem and its weak formulation. Section 3 classifies the interface transmission conditions considered in this work. Section 4 develops the proposed algebraic representation, derives the corresponding global block-matrix system, and discusses its computational properties. Section 5 presents numerical experiments and verification studies. Finally, Section 6 summarizes the main findings and conclusions.

%% file: modelProb.tex
\section{Governing Equations and Model Problem}
\label{sec:interface_classification}

Let $\Omega\subset\mathbb{R}^d$ be a bounded domain partitioned into two non-overlapping subdomains $\Omega^{(1)}$ and $\Omega^{(2)}$, separated by an internal interface $\Gamma$, such that:
\begin{equation} 
\overline{\Omega} = \overline{\Omega^{(1)}}\cup\overline{\Omega^{(2)}},
\qquad \Omega^{(1)}\cap\Omega^{(2)}=\emptyset.
\end{equation}
The external boundary is partitioned into disjoint Dirichlet and Neumann portions, denoted by $\partial\Omega_D$ and $\partial\Omega_N$, respectively. For simplicity, we assume that $\partial\Omega_N\subseteq\partial\Omega^{(2)}$.
We consider a steady diffusion problem governed by:
\begin{equation}
\label{eq:strong_interface_problem}
\begin{aligned}
\nabla\cdot\mathbf q^{(i)} &= f^{(i)},\\
\mathbf q^{(i)} &= -\mathbf k^{(i)}\nabla u^{(i)},
\end{aligned}
\qquad \text{in }\Omega^{(i)},\quad i=1,2,
\end{equation}
where $u^{(i)}$ and $\mathbf q^{(i)}$ denote the solution field and the corresponding flux in subdomain $\Omega^{(i)}$, respectively. Here, $f^{(i)}$ is a prescribed source term and $\mathbf k^{(i)}$ is a symmetric positive-definite diffusion tensor, which may differ between the two subdomains.

\subsection{Boundary conditions}
On the external boundary $\partial\Omega$, we impose homogeneous Dirichlet and prescribed Neumann conditions:
\begin{equation}
\label{eq:boundary}
\begin{cases}
u=0 & \text{on }\partial\Omega_D,\\
\mathbf q^{(2)}\cdot\mathbf n_{\partial\Omega}
=\bar q_0 & \text{on }\partial\Omega_N,
\end{cases}
\end{equation}
where $\mathbf n_{\partial\Omega}$ denotes the outward unit normal
to the external boundary and $\bar q_0$ is the prescribed normal
flux.

\subsection{Interface conditions}
The interface $\Gamma$ separates two regions with potentially different material properties and solution fields. The unit normal $\mathbf n$ to $\Gamma$ is oriented from $\Omega^{(1)}$ towards $\Omega^{(2)}$. Accordingly, the normal flux on either side of the interface is defined with respect to this common orientation as:
\begin{equation}
q_n^{(i)}:=\mathbf q^{(i)}\cdot\mathbf n,
\qquad i=1,2.
\end{equation}
For any piecewise-defined quantity $w$, its jump across the interface is defined by
\begin{equation}
\label{eq:jump_definition}
\llbracket w\rrbracket
:=w^{(2)}-w^{(1)}.
\end{equation}
In particular, the jumps in the solution and normal flux are given by:
\begin{equation}
\label{eq:solution_flux_jumps}
\llbracket u\rrbracket=u^{(2)}-u^{(1)},
\qquad
\llbracket q_n\rrbracket=q_n^{(2)}-q_n^{(1)}.
\end{equation}
The interface behaviour is characterized by the relationship between these two jumps. Depending on the underlying physical model, the solution and flux jumps may be prescribed, constrained to vanish, or related through an interface constitutive law. No particular interface condition is imposed at this stage, allowing the formulation to accommodate different classes of interface behaviour.

 \subsection{The variational formulation}
To derive a general variational formulation, we introduce the broken Sobolev space:
\begin{equation}
V= \left\{ v=(v^{(1)},v^{(2)}):
v^{(i)}\in H^1(\Omega^{(i)}),\
v^{(i)}=0\text{ on }
\partial\Omega_D\cap\partial\Omega^{(i)},~i=1,2
\right\}.
\end{equation}
The functions in this space are defined independently on the two subdomains, and their traces across $\Gamma$ are not constrained to be continuous. Thus, the space accommodates both continuous and discontinuous solution fields, with the specific interface behaviour subsequently determined by the imposed interface conditions.
Multiplying the governing equation in each subdomain by an admissible test function $v^{(i)}$, integrating over $\Omega^{(i)}$, and applying integration by parts yield:
\begin{equation}
\sum_{i=1}^{2} \int_{\Omega^{(i)}}
\mathbf k^{(i)}\nabla u^{(i)}\cdot\nabla v^{(i)}
\,\mathrm d\Omega
+
\int_{\Gamma}
\left(
q_n^{(2)}v^{(2)}-q_n^{(1)}v^{(1)} \right)\,\mathrm d\Gamma = F(v),\label{eq:interface_bilinear_form}
\end{equation}
where the bilinear form and the external loading functional are defined as:
\begin{equation}
\begin{aligned}
a(u,v)
&:=
\sum_{i=1}^{2}
\int_{\Omega^{(i)}}
\mathbf k^{(i)}\nabla u^{(i)}
\cdot\nabla v^{(i)}\,\mathrm d\Omega,
\label{eq:interface_load}\\
F(v) &:= \sum_{i=1}^{2}
\int_{\Omega^{(i)}}f^{(i)}v^{(i)}\,\mathrm d\Omega -
\int_{\partial\Omega_N} \bar q_0v^{(2)}\,\mathrm dS.
\end{aligned}
\end{equation}
The resulting general variational formulation can therefore be written compactly as:
\begin{equation}
\label{eq:general_weak_interface}
a(u,v) + \int_{\Gamma} \left( q_n^{(2)}v^{(2)}-q_n^{(1)}v^{(1)}
\right)\,\mathrm d\Gamma = F(v).
\end{equation}
Importantly, the interface contribution is retained explicitly, rather than being eliminated through an assumed continuity condition. This provides a common starting point for incorporating different interface conditions within the same variational framework.

\begin{rmk}
The interface contribution in \Cref{eq:general_weak_interface} can be expressed in terms of the solution and flux jumps using the algebraic identity:
\begin{equation}
\label{eq:interface_identity}
q_n^{(2)}v^{(2)}-q_n^{(1)}v^{(1)}
=
v^{(2)}\llbracket q_n\rrbracket
+
\llbracket v\rrbracket q_n^{(1)}.
\end{equation}
Consequently, the interface integral admits the equivalent representation
\begin{equation}
\int_{\Gamma}
\left(
q_n^{(2)}v^{(2)}-q_n^{(1)}v^{(1)}
\right)\,\mathrm d\Gamma
=
\int_{\Gamma}
v^{(2)}\llbracket q_n\rrbracket\,\mathrm d\Gamma
+
\int_{\Gamma}
\llbracket v\rrbracket q_n^{(1)}\,\mathrm d\Gamma.
\end{equation}
This decomposition makes explicit the contributions associated with the solution and flux jumps and provides a convenient basis for introducing the different interface conditions considered in the subsequent analysis.
\end{rmk}

%% file: unifiedReform.tex
\section{Classification of Interface Conditions}
Interface conditions can be characterized in terms of the jumps in the primary solution field and its associated normal flux. Various classes of perfect and imperfect interfaces, distinguished by the admissible jumps in these quantities, have been discussed in the literature; see, for instance, Javili et al.~\cite{javilikaessmair2014}. Building upon these established concepts, we organize the interface conditions considered in the present work into seven cases, according to whether the solution and flux jumps are prescribed or remain unknown. The resulting classification is summarized in ~\Cref{tab:interface_conditions}.
\begin{table}[H]
\centering
\setlength{\tabcolsep}{6pt}
\renewcommand{\arraystretch}{1.2}
\begin{tabular}{c|c|c|l|c}
\hline
\textit{Case} &
$\llbracket u\rrbracket$ &
$\llbracket q_n\rrbracket$ &
\textit{Description} &
\textit{Type}\\
\hline\hline
1 & $0$ & $0$ & Perfect interface & Weak\\
\hline
2 & $0$ & $g_q$ & Interface source/sink & Weak\\
\hline
3 & $g_u$ & $0$ & Prescribed solution jump & Strong\\
\hline
4 & $g_u$ & $g_q$ & General prescribed jumps & Strong\\
\hline
5 & Unknown & $0$ & Imperfect interface & Strong\\
\hline
6 & Unknown & $q_n^{(1)}=q_n^{(2)}=0$ &
Impermeable interface & Strong\\
\hline
7 & Unknown & Not necessarily zero &
General interface law & Model dependent\\
\hline
\end{tabular}
\caption{Organization of interface conditions according to the
prescribed or unknown solution and flux jumps. The classification
is based on established interface concepts, see
Javili et al.~\cite{javilikaessmair2014}. Here, $g_u$ and $g_q$ denote
prescribed nonzero functions.}
\label{tab:interface_conditions}
\end{table}

The numerical treatment of interface problems depends on whether the primary solution field is continuous or discontinuous across the
interface. Accordingly, the last column of~\Cref{tab:interface_conditions} distinguishes between weak, strong, and model-dependent discontinuities, based on the regularity of the primary solution field $u$. For weak discontinuities, the solution remains continuous across the interface, and hence a single common set of interface degrees of freedom is sufficient to represent the solution. In contrast, strong discontinuities involve a jump in the solution, requiring independent degrees of freedom to represent the solution traces on either side of the interface. When the solution jump is not explicitly prescribed, as in certain interface constitutive laws, the nature of the discontinuity is determined by the underlying physical model and is therefore classified as model dependent.

In conventional finite element implementations, strong discontinuities are typically represented either through mesh modifications, such as duplicating interface nodes, or through enrichment techniques, such as the extended finite element method (XFEM), which introduce additional approximation functions to capture discontinuous solution fields without requiring mesh conformity to the interface. While enrichment based approaches offer considerable flexibility, their implementation generally involves selecting appropriate enrichment functions and managing the associated enriched degrees of freedom, which may vary depending on the nature of the interface discontinuity and the underlying physical model. To address these issues, the next subsection introduces a topology preserving discretization in which arbitrary solution jumps are represented using auxiliary nodal jump variables, without requiring interface-node duplication or problem-specific enrichment functions.

\section{Unified Finite Element Approximation Across Interfaces}
To capture strong discontinuities without duplicating interface nodes or altering the underlying mesh topology, we introduce a unified jump-enriched finite element approximation. The computational domain $\Omega$ is partitioned by the interface $\Gamma$ into two non-overlapping subdomains, $\Omega_1$ and $\Omega_2$, such that $\Gamma = \partial\Omega_1 \cap \partial\Omega_2$. Instead of assigning two independent sets of nodal unknowns along $\Gamma$, we retain a single primary nodal unknown $u_b^{(2)}$ on the interface and represent the discontinuity explicitly through an auxiliary nodal jump variable $\delta_b = u_b^{(2)} - u_b^{(1)}$. Accordingly, the discrete solution fields in the respective subdomains are approximated as
\begin{equation}
\begin{aligned}
u^{(2)}(\mathbf{x}) &= \sum_{a \in \Omega_2\setminus \Gamma} N_a(\mathbf{x}) u^{(2)}_a + \sum_{b \in \Gamma} N_b(\mathbf{x}) u^{(2)}_b, \qquad \mathbf{x} \in \Omega_2, \\
u^{(1)}(\mathbf{x}) &= \sum_{\alpha \in \Omega_1\setminus\Gamma} N_\alpha(\mathbf{x}) u^{(1)}_\alpha + \sum_{b \in \Gamma} N_b(\mathbf{x}) u_b^{(2)} - \sum_{b\in \Gamma} N_b(\mathbf{x})\delta_b, \qquad \mathbf{x} \in \Omega_1,
\end{aligned}
\label{trail-func}
\end{equation}
where $N_i(\mathbf{x})$ denotes the standard finite element shape functions. Along the interface $\Gamma$, this construction automatically yields the discrete jump, $$\llbracket u \rrbracket = u^{(2)} - u^{(1)} = \sum_{b \in \Gamma} N_b \delta_b.$$ A key advantage of this formulation is that all seven interface cases in Table~\ref{tab:interface_conditions} can be treated within a single algebraic framework: for weak discontinuities (Cases~1--2), $\delta_b$ vanishes identically; for prescribed strong discontinuities (Cases~3--4), $\delta_b$ is enforced directly as a Dirichlet-type constraint; and for unknown jumps (Cases~5--7), $\delta_b$ is solved directly from the coupled weak form, after which the full trace $u^{(1)}$ is recovered through simple post-processing.

\subsection{Unified Weak form and Matrix Assembly}
Consistent with the Bubnov--Galerkin framework, the test functions $v^{(1)}$ and $v^{(2)}$ adopt the exact same jump-enriched expansion as the trial functions:
\begin{equation}
\begin{aligned}
v^{(2)}(\mathbf{x}) &= \sum_{a \in \Omega_2\setminus \Gamma} N_a(\mathbf{x}) v^{(2)}_a + \sum_{b \in \Gamma} N_b(\mathbf{x}) v^{(2)}_b, \qquad \mathbf{x} \in \Omega_2, \\
v^{(1)}(\mathbf{x}) &= \sum_{\alpha \in \Omega_1\setminus\Gamma} N_\alpha(\mathbf{x}) v^{(1)}_\alpha + \sum_{b \in \Gamma} N_b(\mathbf{x}) v_b^{(2)} - \sum_{b\in \Gamma} N_b(\mathbf{x}) \tilde{\delta}_b, \qquad \mathbf{x} \in \Omega_1,
\end{aligned}
\label{test-func}
\end{equation}
so that $\llbracket v \rrbracket = \sum_{b \in \Gamma} N_b \tilde{\delta}_b$ on $\Gamma$. Combining \Cref{eq:interface_bilinear_form,eq:general_weak_interface}, the interface-coupled weak form is :
\begin{equation}
\sum_{i=1}^{2} \int_{\Omega_i} \nabla v^{(i)} \cdot \textbf{k}^{(i)} \nabla u^{(i)} \, \mathrm{d}\Omega +\int_{\Gamma} v^{(2)} \llbracket q_n \rrbracket \, \mathrm{d}\Gamma + \int_{\Gamma} \llbracket v \rrbracket q_n^{(1)} \, \mathrm{d}\Gamma
= \sum_{i=1}^{2} \int_{\Omega_i} v^{(i)} f^{(i)} \, \mathrm{d}\Omega + \int_{\partial\Omega_N} v \, \bar{q}_0 \, \mathrm{d}S.
\label{eq:weak_form_generaluni-1}
\end{equation}
Assume the linear interfacial constitutive law $q_n^{(1)} = q_n^{(2)} = \beta \llbracket u \rrbracket$, then the normal flux jump vanishes identically, that is,~$\llbracket q_n \rrbracket = 0$. Therefore,  the formulation \eqref{eq:weak_form_generaluni-1} simplifies to, 
\begin{equation}
\sum_{i=1}^{2} \int_{\Omega_i} \nabla v^{(i)} \cdot \textbf{k}^{(i)} \nabla u^{(i)} \, \mathrm{d}\Omega  + \int_{\Gamma} \llbracket v \rrbracket  \,\beta \llbracket u \rrbracket \mathrm{d}\Gamma
= \sum_{i=1}^{2} \int_{\Omega_i} v^{(i)} f^{(i)} \, \mathrm{d}\Omega + \int_{\partial\Omega_N} v \, \bar{q}_0 \, \mathrm{d}S.
\label{eq:weak_form_generaluni}
\end{equation}

Now we proceed to describe a  unified block matrix system described as follows.
Let the global vectors of discrete trial and test degrees of freedom be partitioned into interior nodes in $\Omega_1\setminus\Gamma$, interior nodes in $\Omega_2\setminus\Gamma$, shared interface nodes on $\Gamma$, and interfacial jump degrees of freedom on $\Gamma$:
\begin{equation}
\mathbf{d} = 
\begin{bmatrix} 
\mathbf{u}^{(1)}_{\Omega_1} \\ 
\mathbf{u}^{(2)}_{\Omega_2} \\ 
\mathbf{u}^{(2)}_{\Gamma} \\ 
\boldsymbol{\delta}_{\Gamma} 
\end{bmatrix}, 
\qquad 
\mathbf{v} = 
\begin{bmatrix} 
\mathbf{v}^{(1)}_{\Omega_1} \\ 
\mathbf{v}^{(2)}_{\Omega_2} \\ 
\mathbf{v}^{(2)}_{\Gamma} \\ 
\tilde{\boldsymbol{\delta}}_{\Gamma} 
\end{bmatrix}.
\label{eq:global_dof_vectors}
\end{equation}
Since the nodal degrees of freedom on the $\Omega_1$-side of the interface are parametrized via $u^{(1)}_b = u^{(2)}_b - \delta_b$ (with corresponding test variations $v^{(1)}_b = v^{(2)}_b - \tilde{\delta}_b$), substituting the discrete trial 
 \eqref{trail-func}~and test expansions \eqref{test-func} into \eqref{eq:weak_form_generaluni} yields the symmetric global linear system:
\begin{equation}
\begin{bmatrix}
\mathbf{K}^{(1)}_{\Omega_1 \Omega_1} & \mathbf{0} & \mathbf{K}^{(1)}_{\Omega_1 \Gamma} & -\mathbf{K}^{(1)}_{\Omega_1 \Gamma} \\[4pt]
\mathbf{0} & \mathbf{K}^{(2)}_{\Omega_2 \Omega_2} & \mathbf{K}^{(2)}_{\Omega_2 \Gamma} & \mathbf{0} \\[4pt]
\mathbf{K}^{(1)}_{\Gamma \Omega_1} & \mathbf{K}^{(2)}_{\Gamma \Omega_2} & \mathbf{K}^{(1)}_{\Gamma \Gamma} + \mathbf{K}^{(2)}_{\Gamma \Gamma} & -\mathbf{K}^{(1)}_{\Gamma \Gamma} \\[4pt]
-\mathbf{K}^{(1)}_{\Gamma \Omega_1} & \mathbf{0} & -\mathbf{K}^{(1)}_{\Gamma \Gamma} & \mathbf{K}^{(1)}_{\Gamma \Gamma} + \mathbf{M}^{\beta}_{\Gamma\Gamma}
\end{bmatrix}
\begin{bmatrix}
\mathbf{u}^{(1)}_{\Omega_1} \\[4pt]
\mathbf{u}^{(2)}_{\Omega_2} \\[4pt]
\mathbf{u}^{(2)}_{\Gamma} \\[4pt]
\boldsymbol{\delta}_{\Gamma}
\end{bmatrix}
=
\begin{bmatrix}
\mathbf{f}^{(1)}_{\Omega_1} \\[4pt]
\mathbf{f}^{(2)}_{\Omega_2} \\[4pt]
\mathbf{f}^{(1)}_{\Gamma} + \mathbf{f}^{(2)}_{\Gamma} \\[4pt]
-\mathbf{f}^{(1)}_{\Gamma}
\end{bmatrix},
\label{eq:block_system}
\end{equation}
where the subdomain stiffness submatrices, interfacial conductance matrix, and consistent load subvectors are defined as:
\begin{equation}
\begin{aligned}
\mathbf{K}^{(i)}_{IJ} &= \int_{\Omega_i} \nabla \mathbf{N}_I^{\mathsf{T}} \, \mathbf{k}^{(i)} \nabla \mathbf{N}_J \, \mathrm{d}\Omega, \qquad (i = 1, 2),\\[4pt] %\label{eq:stiffness_submatrices} \\[4pt]
\mathbf{M}^{\beta}_{\Gamma\Gamma} &= \int_{\Gamma} \beta \, \mathbf{N}_{\Gamma}^{\mathsf{T}} \mathbf{N}_{\Gamma} \, \mathrm{d}\Gamma, \\[4pt]%\label{eq:interface_mass_matrix} \\[4pt]
\mathbf{f}^{(i)}_{I} &= \int_{\Omega_i} \mathbf{N}_I^{\mathsf{T}} f^{(i)} \, \mathrm{d}\Omega + \int_{\partial\Omega_N \cap \partial\Omega_i} \mathbf{N}_I^{\mathsf{T}} \bar{q}_0 \, \mathrm{d}S, %\label{eq:load_subvectors}
\end{aligned}
\label{eqn:stiffness_load_def}
\end{equation}
with index sets $I \in \{\Omega_i\setminus\Gamma~i=1,2\}$~and~$J \in \{ \Gamma\}$ denoting the interior and interface shape function row arrays, respectively.
 \subsection{Unified Treatment of Interface Cases and Post-Processing}

The block structure of \Cref{eq:block_system} accommodates all seven cases in Table~\ref{tab:interface_conditions} without re-meshing or modifying the global connectivity:
\begin{itemize}
    \item \textbf{Weak discontinuities (Cases~1 and~2):} The jump vector is constrained to $\boldsymbol{\delta} = \mathbf{0}$, eliminating the fourth block row and column and reducing \Cref{eq:block_system} to a standard continuous finite element system with $\mathbf{g}_q = \mathbf{0}$ (Case~1) or $\mathbf{g}_q \neq \mathbf{0}$ (Case~2).
    \item \textbf{Prescribed strong discontinuities (Cases~3 and~4):} The jump vector is prescribed directly as $\boldsymbol{\delta} = \mathbf{g}_u$ on $\Gamma$ and moved to the right-hand side as a standard Dirichlet condition.
    \item \textbf{Imperfect and impermeable interfaces (Cases~5 and~6):} Solving \Cref{eq:block_system} directly yields both $\mathbf{u}_{\Gamma}^{(2)}$ and $\boldsymbol{\delta}$ for a given $\beta > 0$ (Case~5) or $\beta = 0$ (Case~6, where $\mathbf{M}_{\Gamma\Gamma}^{\beta} = \mathbf{0}$).
    \item \textbf{General interface constitutive laws (Case~7):} Any linear or nonlinear relation $q_n^{(1)} = \mathcal{F}(u^{(2)}, \llbracket u \rrbracket)$ is incorporated solely by updating the interface contribution in the fourth block row (see \Cref{sec:nonlinear-jump}).
\end{itemize}

Once the global system in \Cref{eq:block_system} is solved, the nodal solution trace on the $\Omega_1$ side of the interface is recovered in post-processing via the algebraic relation
\begin{equation}
\mathbf{u}_{\Gamma}^{(1)} = \mathbf{u}_{\Gamma}^{(2)} - \boldsymbol{\delta},
\end{equation}
providing complete, discontinuous field representations across both subdomains.

\subsection{Extension to Nonlinear Interfacial Constitutive Laws}
\label{sec:nonlinear-jump}
Consider a general nonlinear imperfect interface characterized by continuous normal flux across $\Gamma$, $q_n^{(1)} = q_n^{(2)} = \mathcal{F}(u^{(2)}, \llbracket u \rrbracket)$, where $\mathcal{F}$ is a smooth, nonlinear function of the interface trace $u^{(2)}|_{\Gamma}$ and the field jump $\llbracket u \rrbracket$ (such as radiative heat transfer, nonlinear contact resistance, or Butler--Volmer interfacial kinetics). Since $\llbracket q_n \rrbracket = 0$ still holds identically, substituting the discrete trial and test expansions into \eqref{eq:weak_form_generaluni} yields the nonlinear residual system $\mathbf{R}(\mathbf{u}) = \mathbf{0}$:
\begin{equation}
\mathbf{R}(\mathbf{u}) =
\begin{bmatrix}
\mathbf{K}^{(1)}_{\Omega_1 \Omega_1} \mathbf{u}^{(1)}_{\Omega_1} + \mathbf{K}^{(1)}_{\Omega_1 \Gamma} \mathbf{u}^{(2)}_{\Gamma} - \mathbf{K}^{(1)}_{\Omega_1 \Gamma} \boldsymbol{\delta}_{\Gamma} - \mathbf{f}^{(1)}_{\Omega_1} \\[4pt]
\mathbf{K}^{(2)}_{\Omega_2 \Omega_2} \mathbf{u}^{(2)}_{\Omega_2} + \mathbf{K}^{(2)}_{\Omega_2 \Gamma} \mathbf{u}^{(2)}_{\Gamma} - \mathbf{f}^{(2)}_{\Omega_2} \\[4pt]
\mathbf{K}^{(1)}_{\Gamma \Omega_1} \mathbf{u}^{(1)}_{\Omega_1} + \mathbf{K}^{(2)}_{\Gamma \Omega_2} \mathbf{u}^{(2)}_{\Omega_2} + \left(\mathbf{K}^{(1)}_{\Gamma \Gamma} + \mathbf{K}^{(2)}_{\Gamma \Gamma}\right) \mathbf{u}^{(2)}_{\Gamma} - \mathbf{K}^{(1)}_{\Gamma \Gamma} \boldsymbol{\delta}_{\Gamma} - \left(\mathbf{f}^{(1)}_{\Gamma} + \mathbf{f}^{(2)}_{\Gamma}\right) \\[4pt]
-\mathbf{K}^{(1)}_{\Gamma \Omega_1} \mathbf{u}^{(1)}_{\Omega_1} - \mathbf{K}^{(1)}_{\Gamma \Gamma} \mathbf{u}^{(2)}_{\Gamma} + \mathbf{K}^{(1)}_{\Gamma \Gamma} \boldsymbol{\delta}_{\Gamma} + \mathbf{r}^{\mathrm{int}}_{\Gamma}(\mathbf{u}^{(2)}_{\Gamma}, \boldsymbol{\delta}_{\Gamma}) + \mathbf{f}^{(1)}_{\Gamma}
\end{bmatrix}
= \mathbf{0},
\label{eq:nonlinear_residual_system}
\end{equation}
where the nonlinear interfacial flux vector is defined as
\begin{equation}
\mathbf{r}^{\mathrm{int}}_{\Gamma}(\mathbf{u}^{(2)}_{\Gamma}, \boldsymbol{\delta}_{\Gamma}) = \int_{\Gamma} \mathbf{N}_{\Gamma}^{\mathsf{T}} \, \mathcal{F}\!\left(\mathbf{N}_{\Gamma}\mathbf{u}^{(2)}_{\Gamma}, \, \mathbf{N}_{\Gamma}\boldsymbol{\delta}_{\Gamma}\right) \mathrm{d}\Gamma.
\label{eq:nonlinear_interface_vector}
\end{equation}

\paragraph{Consistent Newton--Raphson Linearization}
Solving \eqref{eq:nonlinear_residual_system} via the Newton--Raphson method requires the iterative update $^{(k+1)}\mathbf{d} = ^{(k)}\mathbf{d} + \Delta ^{(k)}\mathbf{d}$, where the incremental vector $\Delta ^{(k)}\mathbf{d}$ satisfies the linearized block system $\mathbf{J}(^{(k)}\mathbf{d})\,\Delta^{(k)}\mathbf{d} = -\mathbf{R}( ^{(k)}\mathbf{d})$:
\begin{equation}
\begin{bmatrix}
\mathbf{K}^{(1)}_{\Omega_1 \Omega_1} & \mathbf{0} & \mathbf{K}^{(1)}_{\Omega_1 \Gamma} & -\mathbf{K}^{(1)}_{\Omega_1 \Gamma} \\[4pt]
\mathbf{0} & \mathbf{K}^{(2)}_{\Omega_2 \Omega_2} & \mathbf{K}^{(2)}_{\Omega_2 \Gamma} & \mathbf{0} \\[4pt]
\mathbf{K}^{(1)}_{\Gamma \Omega_1} & \mathbf{K}^{(2)}_{\Gamma \Omega_2} & \mathbf{K}^{(1)}_{\Gamma \Gamma} + \mathbf{K}^{(2)}_{\Gamma \Gamma} & -\mathbf{K}^{(1)}_{\Gamma \Gamma} \\[4pt]
-\mathbf{K}^{(1)}_{\Gamma \Omega_1} & \mathbf{0} & -\mathbf{K}^{(1)}_{\Gamma \Gamma} + \mathbf{J}_{\delta u} & \mathbf{K}^{(1)}_{\Gamma \Gamma} + \mathbf{J}_{\delta \delta}
\end{bmatrix}
\begin{bmatrix}
\Delta\mathbf{u}^{(1)}_{\Omega_1} \\[4pt]
\Delta\mathbf{u}^{(2)}_{\Omega_2} \\[4pt]
\Delta\mathbf{u}^{(2)}_{\Gamma} \\[4pt]
\Delta\boldsymbol{\delta}_{\Gamma}
\end{bmatrix}
= -\mathbf{R}(\mathbf{d}^{(k)}).
\label{eq:nonlinear_jacobian_system}
\end{equation}
Here, the consistent interfacial tangent sub-matrices are given by:
 \begin{equation}
\mathbf{J}_{\delta u} = \frac{\partial \mathbf{r}^{\mathrm{int}}_{\Gamma}}{\partial \mathbf{u}^{(2)}_{\Gamma}} \qquad\text{and}\qquad
\mathbf{J}_{\delta \delta} = \frac{\partial \mathbf{r}^{\mathrm{int}}_{\Gamma}}{\partial \boldsymbol{\delta}_{\Gamma}}.
 \end{equation}
% \begin{subequations}
% \begin{align}
% \mathbf{J}_{\delta u} &= \frac{\partial \mathbf{r}^{\mathrm{int}}_{\Gamma}}{\partial \mathbf{u}^{(2)}_{\Gamma}} = \int_{\Gamma} \frac{\partial \mathcal{F}}{\partial u^{(2)}} \, \mathbf{N}_{\Gamma}^{\mathsf{T}} \mathbf{N}_{\Gamma} \, \mathrm{d}\Gamma, \label{eq:jacobian_du} \\[4pt]
% \mathbf{J}_{\delta \delta} &= \frac{\partial \mathbf{r}^{\mathrm{int}}_{\Gamma}}{\partial \boldsymbol{\delta}_{\Gamma}} = \int_{\Gamma} \frac{\partial \mathcal{F}}{\partial \delta} \, \mathbf{N}_{\Gamma}^{\mathsf{T}} \mathbf{N}_{\Gamma} \, \mathrm{d}\Gamma. \label{eq:jacobian_dd}
% \end{align}
% \end{subequations}
Notice that the bulk stiffness blocks remain constant throughout the iterations, while only the two interface tangent blocks $\mathbf{J}_{\delta u}$ and $\mathbf{J}_{\delta \delta}$ require reassembly along $\Gamma$. Furthermore, whenever $\partial \mathcal{F}/\partial u^{(2)} \neq 0$, the term $\mathbf{J}_{\delta u}$ breaks the symmetry of the global Jacobian matrix $\mathbf{J}$.

 \subsection{Computational advantages of the proposed formulation}

The jump-parametrized representation offers several computational advantages
arising from the explicit separation of the bulk field and the interfacial
jump. In particular, the following features are noteworthy:
\begin{enumerate}
    \item \emph{Explicit representation of the interface jump:}
    The jump variable \(\boldsymbol{\delta}_{\Gamma}\) provides direct
    access to the interfacial discontinuity and allows the constitutive
    response to be incorporated explicitly into the algebraic system.
    This is particularly useful when the interface law depends nonlinearly
    on the temperature jump.

    \item \emph{Treatment of different interface conditions:}
    Prescribed jumps and constitutive transmission conditions can be
    incorporated by modifying the interface-associated equations, without
    changing the underlying bulk discretization. In particular, perfect
    contact can be imposed strongly through
    \(\boldsymbol{\delta}_{\Gamma}=\mathbf{0}\), avoiding the need for
    large penalty parameters.

    \item \emph{Primal algebraic structure:}
    The interface constitutive relation is incorporated directly into
    the primal system without introducing independent flux multipliers.
    Consequently, the formulation does not require a separate
    multiplier space or the associated inf-sup stability analysis.

    \item \emph{Modular nonlinear assembly:}
    For nonlinear interface laws, the nonlinear residual and its
    consistent tangent are assembled through the interface-associated
    contributions. For linear bulk problems, the bulk stiffness
    contributions remain unchanged during Newton iterations, thereby
    localizing the nonlinear updates to the interface.

    \item \emph{Compatibility with different discretizations:}
    The jump decomposition is independent of the particular bulk
    approximation and can be incorporated into compatible finite element
    formulations, including FEM, IGA, polygonal finite element methods
    and VEM.
\end{enumerate}

%% file: onedAnal.tex
\section{One-dimensional heat conduction across an interface: analytical solution}
\label{appx:oned_anal}

Consider the one-dimensional steady-state heat conduction problem
\begin{equation}
-\frac{\mathrm{d}}{\mathrm{d}x}
\left(k(x)\frac{\mathrm{d}u}{\mathrm{d}x}\right)=0,
\qquad x\in(0,1)\setminus\{x_b\},
\label{eqn:oned_ode}
\end{equation}
where \(x_b\) denotes the location of the internal interface separating two
materials with piecewise constant thermal conductivities:
\begin{equation}
k(x)=
\begin{cases}
k^{(1)}, & 0<x<x_b,\\
k^{(2)}, & x_b<x<1.
\end{cases}
\end{equation}
The prescribed boundary temperatures are
\begin{equation}
u^{(1)}(0)=u_L,
\qquad
u^{(2)}(1)=u_R.
\label{eqn:oned_bcs}
\end{equation}
The heat flux in each subdomain is defined by
\begin{equation}
q^{(i)}=-k^{(i)}\frac{\mathrm{d}u^{(i)}}{\mathrm{d}x},
\qquad i=1,2,
\end{equation}
with the positive direction taken along the increasing \(x\)-axis. Accordingly,
the interface jumps are defined as
\begin{equation}
\llbracket u\rrbracket
=u^{(2)}(x_b)-u^{(1)}(x_b),
\qquad
\llbracket q\rrbracket
=q^{(2)}-q^{(1)}.
\end{equation}
For convenience, introduce the subdomain lengths, thermal resistances, and
auxiliary quantities
\begin{equation}
\begin{aligned}
h_1 &= x_b, 
& h_2 &= 1-x_b,\\
r_1 &= \frac{h_1}{k^{(1)}},
& r_2 &= \frac{h_2}{k^{(2)}},\\
S &= r_1+r_2,
& D &= u_R-u_L.
\end{aligned}
\label{eqn:oned_auxiliary}
\end{equation}

\subsection*{Prescribed temperature and heat-flux jumps}

Let the prescribed temperature and heat-flux jumps be
\begin{equation}
\llbracket u\rrbracket=g_u,
\qquad
\llbracket q\rrbracket=g_q.
\end{equation}
Since the heat flux is constant within each subdomain, integration of
\eqref{eqn:oned_ode}, together with the boundary conditions
\eqref{eqn:oned_bcs}, gives the temperature distribution
\begin{equation}
u(x)=
\begin{cases}
u^{(1)}(x)
=u_L-\dfrac{x}{k^{(1)}}q^{(1)},
& 0\leq x<x_b,\\[1em]
u^{(2)}(x)
=u_R+\dfrac{1-x}{k^{(2)}}q^{(2)},
& x_b<x\leq 1.
\end{cases}
\label{eqn:oned_temperature_general}
\end{equation}
The prescribed flux-jump condition yields
\begin{equation}
q^{(2)}=q^{(1)}+g_q.
\label{eqn:oned_flux_jump}
\end{equation}
Evaluating the temperature jump at \(x=x_b\) and using
\eqref{eqn:oned_temperature_general}, we obtain
\begin{equation}
\begin{aligned}
g_u
&=u^{(2)}(x_b)-u^{(1)}(x_b)\\
&=D+r_1q^{(1)}+r_2q^{(2)}\\
&=D+Sq^{(1)}+r_2g_q.
\end{aligned}
\end{equation}
Consequently, the two constant heat fluxes are given by
\begin{equation} 
q^{(1)}
=\frac{g_u-D-r_2g_q}{S}\qquad\text{and}\qquad
q^{(2)}
=q^{(1)}+g_q. 
\label{eqn:oned_flux_general}
\end{equation}
Substitution into \eqref{eqn:oned_temperature_general} provides the analytical
temperature field for arbitrary prescribed temperature and heat-flux jumps.

% The four prescribed-jump interface conditions considered in this work follow
% as particular cases: ~(i)~\textit{Perfect thermal contact:}
%     \(g_u=0\) and \(g_q=0\);~(ii)~\textit{Prescribed heat-flux jump:}
%     \(g_u=0\) and \(g_q\neq 0\);~(iii)~
%  \textit{Prescribed temperature jump:}
%     \(g_u\neq 0\) and \(g_q=0\);~(iv)~\textit{Prescribed temperature and heat-flux jumps:}both \(g_u\) and \(g_q\) are nonzero.
\subsection*{Kapitza interface condition}
For a Kapitza interface, the heat flux is continuous across the interface,
whereas the temperature may be discontinuous. Thus,
\begin{equation}
q^{(1)}=q^{(2)}=q,
\qquad
\llbracket u\rrbracket=-\beta q,
\label{eqn:oned_kapitza_law}
\end{equation}
where \(\beta>0\) denotes the interfacial thermal resistance. Integrating
\eqref{eqn:oned_ode} in the two subdomains gives
\begin{equation}
\begin{aligned}
u^{(1)}(x)
=u_L-\frac{q}{k^{(1)}}x\qquad\text{and}\qquad
u^{(2)}(x)
=u_R+\frac{q}{k^{(2)}}(1-x).
\end{aligned}
\label{eqn:oned_kapitza_temperature}
\end{equation}
The resulting temperature jump at the interface is
\begin{equation}
\begin{aligned}
\llbracket u\rrbracket
&=u^{(2)}(x_b)-u^{(1)}(x_b)\\
&=(u_R-u_L)+q(r_1+r_2)\,=\,D+Sq.
\end{aligned}
\label{eqn:oned_kapitza_jump}
\end{equation}
Combining \eqref{eqn:oned_kapitza_law} and
\eqref{eqn:oned_kapitza_jump}, we obtain
\begin{equation}
D+Sq=-\beta q.
\end{equation}
Therefore, the heat flux and the corresponding temperature jump are
\begin{equation} 
q=-\frac{D}{S+\beta} \qquad\text{and}\qquad \llbracket u\rrbracket=\frac{\beta D}{S+\beta}. 
\label{eqn:oned_kapitza_solution}
\end{equation}
In the limit \(\beta\to0\), the temperature jump vanishes and the solution
reduces to that of perfect thermal contact.